\documentclass[preprint,11pt]{elsarticle}
\usepackage[bookmarks=true,colorlinks=true,linkcolor=blue]{hyperref}

\biboptions{sort&compress}

\usepackage{calc}
\usepackage[margin=1.in]{geometry}

\usepackage[usenames,dvipsnames,svgnames,table]{xcolor}

\usepackage{graphicx}

\usepackage{amsmath,amssymb}
\usepackage{verbatim}
\usepackage{float}
\usepackage{caption} 
\usepackage{tikz}
\usepackage{hyperref}

\newlength{\tfwidth}
\newlength{\tfheight}
\newlength{\tfxa}
\newlength{\tfxb}
\newlength{\tfya}
\newlength{\tfyb}
\newcommand{\trimFigWithBox}[6]{%
\setlength\fboxsep{0pt}%
\setlength\fboxrule{1.0pt}
\fbox{\includegraphics[width=#2, clip, trim=#3 #4 #5 #6]{#1}}%
}
\newcommand{\trimFigNoBox}[6]{%
\setlength\fboxsep{1pt}
\setlength\fboxrule{0.0pt}
\fbox{\includegraphics[width=#2, clip, trim=#3 #4 #5 #6]{#1}}%
}
\newcommand{\trimFigHeightWithBox}[6]{%
\setlength\fboxsep{0pt}%
\setlength\fboxrule{1.0pt}
\fbox{\includegraphics[height=#2, clip, trim=#3 #4 #5 #6]{#1}}%
}
\newcommand{\trimFigHeightNoBox}[6]{%
\setlength\fboxsep{1pt}
\setlength\fboxrule{0.0pt}
\fbox{\includegraphics[height=#2, clip, trim=#3 #4 #5 #6]{#1}}%
}

\newsavebox\figBox

\newcommand{\trimw}[6]{%
\sbox\figBox{\includegraphics{#1}}
\setlength{\tfwidth}{\the\wd\figBox}
\setlength{\tfheight}{\the\ht\figBox}
\setlength{\tfxa}{\tfwidth*\real{#3}}%
\setlength{\tfxb}{\tfwidth*\real{#4}}%
\setlength{\tfya}{\tfheight*\real{#5}}%
\setlength{\tfyb}{\tfheight*\real{#6}}%
\trimFigNoBox{#1}{#2}{\tfxa}{\tfya}{\tfxb}{\tfyb}%
}

\newcommand{\trimwb}[6]{%

\sbox\figBox{\includegraphics{#1}}
\setlength{\tfwidth}{\the\wd\figBox}
\setlength{\tfheight}{\the\ht\figBox}
\setlength{\tfxa}{\tfwidth*\real{#3}}%
\setlength{\tfxb}{\tfwidth*\real{#4}}%
\setlength{\tfya}{\tfheight*\real{#5}}%
\setlength{\tfyb}{\tfheight*\real{#6}}%
\trimFigWithBox{#1}{#2}{\tfxa}{\tfya}{\tfxb}{\tfyb}%
}

\newcommand{\trimh}[6]{%
\sbox\figBox{\includegraphics{#1}}
\setlength{\tfwidth}{\the\wd\figBox}
\setlength{\tfheight}{\the\ht\figBox}
\setlength{\tfxa}{\tfwidth*\real{#3}}%
\setlength{\tfxb}{\tfwidth*\real{#4}}%
\setlength{\tfya}{\tfheight*\real{#5}}%
\setlength{\tfyb}{\tfheight*\real{#6}}%
\trimFigHeightNoBox{#1}{#2}{\tfxa}{\tfya}{\tfxb}{\tfyb}%
}

\newcommand{\trimhb}[6]{%

\sbox\figBox{\includegraphics{#1}}
\setlength{\tfwidth}{\the\wd\figBox}
\setlength{\tfheight}{\the\ht\figBox}
\setlength{\tfxa}{\tfwidth*\real{#3}}%
\setlength{\tfxb}{\tfwidth*\real{#4}}%
\setlength{\tfya}{\tfheight*\real{#5}}%
\setlength{\tfyb}{\tfheight*\real{#6}}%
\trimFigHeightWithBox{#1}{#2}{\tfxa}{\tfya}{\tfxb}{\tfyb}%
}
\usepackage{xargs}

\newcommandx{\figByHeight}[9][5=0, 6=0, 7=0, 8=0,9=]{
\draw (#1,#2) node[anchor=south west,xshift=-16pt,yshift=-4pt] {\trimh{#3}{#4}{#5}{#6}{#7}{#8}};}
\newcommandx{\figByHeightb}[9][5=0, 6=0, 7=0, 8=0,9=]{
\draw (#1,#2) node[anchor=south west,xshift=-16pt,yshift=-4pt] {\trimhb{#3}{#4}{#5}{#6}{#7}{#8}};}

\newcommandx{\figByHeightWithLabel}[9][5=0, 6=0, 7=0, 8=0,9=]{
\draw (#1,#2) node[anchor=south west,xshift=-16pt,yshift=-4pt] {\trimh{#3}{#4}{#5}{#6}{#7}{#8}} node[draw=white,fill=white,inner sep=1pt,anchor=south west] {#9};}
\newcommandx{\figByHeightWithLabelb}[9][5=0, 6=0, 7=0, 8=0,9=]{
\draw (#1,#2) node[anchor=south west,xshift=-16pt,yshift=-4pt] {\trimhb{#3}{#4}{#5}{#6}{#7}{#8}} node[draw=white,fill=white,inner sep=1pt,anchor=south west] {#9};}

\newcommandx{\figByWidth}[9][5=0, 6=0, 7=0, 8=0,9=]{
\draw (#1,#2) node[anchor=south west,xshift=-16pt,yshift=-4pt] {\trimw{#3}{#4}{#5}{#6}{#7}{#8}};}
\newcommandx{\figByWidthb}[9][5=0, 6=0, 7=0, 8=0,9=]{
\draw (#1,#2) node[anchor=south west,xshift=-16pt,yshift=-4pt] {\trimwb{#3}{#4}{#5}{#6}{#7}{#8}};}

\newcommandx{\figByWidthWithLabel}[9][5=0, 6=0, 7=0, 8=0,9=]{
\draw (#1,#2) node[anchor=south west,xshift=-16pt,yshift=-4pt] {\trimw{#3}{#4}{#5}{#6}{#7}{#8}} node[draw=white,fill=white,inner sep=1pt,anchor=south west] {#9};}
\newcommandx{\figByWidthWithLabelb}[9][5=0, 6=0, 7=0, 8=0,9=]{
\draw (#1,#2) node[anchor=south west,xshift=-16pt,yshift=-4pt] {\trimwb{#3}{#4}{#5}{#6}{#7}{#8}} node[draw=white,fill=white,inner sep=1pt,anchor=south west] {#9};}

\usepackage{algorithm,algpseudocode}
\usepackage{algorithmicx}
\algrenewcommand\alglinenumber[1]{\footnotesize #1:} 

\usepackage{todonotes}

\definecolor{purple}{rgb}{.7, 0., .8}

\newcommand{\blue}{\color{blue}}

\newcommand{\ssf}{\scriptscriptstyle}

\newcommand{\smallss}{\sffamily\small}

\newcommand{\half}{\frac{1}{2}}

\newcommand{\Ih}{\mathcal{I}_{\ssf H}}
\newcommand{\Th}{\mathcal{T}_{\ssf H}}
\newcommand{\Bh}{\mathcal{B}_{\ssf H}}

\newcommand{\dt}{{\Delta t}}
\newcommand{\dr}{{\Delta r}}
\newcommand{\ds}{{\Delta s}}
\newcommand{\dx}{{\Delta x}}
\newcommand{\dy}{{\Delta y}}

\newcommand{\Lh}{L_{\ssf H}}
\newcommand{\Nh}{N_{\ssf H}}

\newtheorem{theorem}{Theorem}

\newenvironment{proof}[1][Proof]{\begin{trivlist}
\item[\hskip \labelsep {\bfseries #1.}]}{\end{trivlist}}

\newcommand{\bni}{\bigskip\noindent}
\newcommand{\mni}{\medskip\noindent}

\newcommand{\p}{\partial}

\newcommand{\f}[2]{\frac{#1}{#2}}

\def\ba#1\ea{\begin{align}#1\end{align}}

\def\bas#1\eas{\begin{align*}#1\end{align*}}

\def\bat#1\eat{\begin{alignat}{3}#1\end{alignat}}

\def\bats#1\eats{\begin{alignat*}{3}#1\end{alignat*}}

\newcommand{\bse}{\begin{subequations}}
\newcommand{\ese}{\end{subequations}}

\newcommand{\eqdef}{\overset{{\rm def}}{=}}

\newcommand{\bv}{\mathbf{ b}}

\newcommand{\dv}{\mathbf{ d}}
\newcommand{\ev}{\mathbf{ e}}
\newcommand{\fv}{\mathbf{ f}}

\newcommand{\iv}{\mathbf{ i}}
\newcommand{\jv}{\mathbf{ j}}
\newcommand{\kv}{\mathbf{ k}}

\newcommand{\nv}{\mathbf{ n}}

\newcommand{\rv}{\mathbf{ r}}

\newcommand{\uv}{\mathbf{ u}}
\newcommand{\vv}{\mathbf{ v}}
\newcommand{\wv}{\mathbf{ w}}
\newcommand{\xv}{\mathbf{ x}}

\newcommand{\zv}{\mathbf{ z}}

\newcommand{\Gv}{\mathbf{ G}}

\newcommand{\zerov}{\mathbf{0}}

\newcommand{\grad}{\nabla}

\newcommand{\Real}{{\mathbb R}}
\newcommand{\Bc}{{\mathcal B}}

\newcommand{\Mc}{{\mathcal M}}
\newcommand{\Nc}{{\mathcal N}}

\newcommand{\dBar}{\bar{d}}

\newcommand{\dvBar}{\bar{\dv}}
\newcommand{\nBar}{\bar{n}}
\newcommand{\wBar}{\bar{w}}

\newcommand{\tcr}{\gamma}

\newcommand{\ts}{\beta}
\newcommand{\cfl}{C_{\ssf CFL}}
\newcommand{\dxMin}{\dx_{\ssf min}}

\usepackage{fancyvrb}
\usepackage{empheq}
\usepackage[most]{tcolorbox}
\newtcbox{\mymath}[1][]{%
nobeforeafter, math upper, tcbox raise base,
enhanced, colframe=blue!60!black,
colback=blue!20, boxrule=1pt,
#1}
\usepackage{mdframed}
\newcommand{\shadedBoxWithShadow}[3]{\begin{empheq}[box={\mymath[colback=#2!20,colframe=#2!60!black,drop lifted shadow, sharp corners]}]{#1} #3\end{empheq}}

\usepackage{listings}
\newcommand{\REV}[2]{%
  \ifthenelse{\equal{#1}{0}}{{{\color{black}#2}}}{}%
  \ifthenelse{\equal{#1}{1}}{{{\color{black}#2}}}{}%
  \ifthenelse{\equal{#1}{2}}{{{\color{black}#2}}}{}%
  \ifthenelse{\equal{#1}{3}}{{{\color{black}#2}}}{}%
  \ifthenelse{\equal{#1}{4}}{{{\color{black}#2}}}{}%
  \ifthenelse{\equal{#1}{5}}{{{\color{black}#2}}}{}%
}

\begin{document}

\begin{frontmatter}
\title{High Order Accurate Hermite Schemes on Curvilinear Grids with Compatibility Boundary Conditions}

\author[lanl]{Allen Alvarez Loya\corref{cor}\fnref{AllenThanks}}
\ead{aalvarezloya@lanl.gov}

\author[vtu]{Daniel Appel\"o\fnref{DanielThanks}}
\ead{appelo@vt.edu}

\author[rpi]{William D.~Henshaw\fnref{BillThanks}}
\ead{henshw@rpi.edu}

\cortext[cor]{Corresponding author}

\address[lanl]{Los Alamos National Laboratory, Los Alamos, NM 87544, USA}
\fntext[AllenThanks]{Research supported by National Science Foundation under grants DGE-1650115 and DMS-2213261}

\address[vtu]{Department of Mathematics, Virginia Tech, Blacksburg, VA 24061 USA}
\fntext[DanielThanks]{Research supported by National Science Foundation under grant DMS-2345225, and Virginia Tech.}

\address[rpi]{Department of Mathematical Sciences, Rensselaer Polytechnic Institute, Troy, NY 12180, USA}
\fntext[BillThanks]{Research supported by the National Science Foundation under grants DMS-1519934 and DMS-1818926.}


\begin{abstract} 
High order accurate Hermite methods for the wave equation on curvilinear domains are presented.
Boundaries are treated using centered compatibility conditions rather than more standard one-sided approximations. 
Both first-order-in-time (FOT) and second-order-in-time (SOT) Hermite schemes are developed.
Hermite methods use the solution and multiple derivatives as unknowns and achieve
space-time orders of accuracy $2m-1$ (FOT) and $2m$ (SOT) for methods using $(m+1)^d$ degree of freedom per node in $d$ dimensions. 
The compatibility boundary conditions (CBCs) are based on taking time derivatives of the boundary conditions
and using the governing equations to 
replace the time derivatives with spatial derivatives.
These resulting constraint equations augment the Hermite scheme on the boundary.
The solvability of the equations resulting from the compatibility conditions are analyzed.
Numerical examples demonstrate the accuracy and stability of the new schemes in two dimensions.
\end{abstract}


\end{frontmatter}

\clearpage
\tableofcontents

\clearpage

\section{Introduction}\label{sec:Introduction}

We develop high-order accurate Hermite methods for the wave equation on curvilinear grids.
Both first-order in time (FOT) and \REV{3}{second-order in time (SOT) schemes are developed}.
Compatibility boundary conditions (CBCs) are used
to give high-order accurate centered approximations to boundary conditions rather than more common one-sided approximations;
these centered boundary conditions are generally more stable and accurate than
using one sided approximations~\cite{hassanieh2021local}.
A key result of the current article is to show how CBCs can be incorporated into Hermite schemes. \REV{5}{This is the first Hermite method that can handle 
complex geometry \textit{and} boundary conditions purely within the Hermite method framework.} The Hermite method approximates the solution to a partial differential equation
using degrees of freedom at each node representing the solution and derivatives up to degree $m$ resulting in 
$(m+1)^d$ degrees of freedom per node in $d$ dimensions.
The nodal values are interpolated to the cell centers and advanced in time using a Taylor series in time over half a time step.  
The cell centered values at the half time-step are then interpolated to the nodes and advanced the second-half time step. The evolution of the degrees of freedom is local, minimizing the communication and storage costs. 
The resulting schemes have order of accuracy $2m-1$ for the FOT scheme, and $2m$ for the \REV{3}{SOT} scheme.
We note that when the FOT scheme is used on Cartesian grids and with $d(2m+1)$ terms in the Taylor series in time (see \cite{goodrich2006hermite,appelo2018hermite})  the schemes have a CFL number that is one, i.e.  they are stable for $c\dt/h \le 1$ where $c$ is the wave speed and $h$ in the grid spacing. This results holds at any order of accuracy and is significantly better than many other high-order schemes. To achieve space-time accuracies of $2m-1$ for the FOT scheme it is sufficient to use a Taylor series in time with $2m+1$ terms, which is what we do here. This reduces the time-step slightly, here we use $c\dt/h = 0.5$,  but does not significantly change the efficiency of the methods as the cost per time-step is reduced by a factor $d$. The \REV{3}{SOT} scheme is even more efficient, requiring only $m+1$ terms in time for any dimension on Cartesian meshes. On curvilinear meshes we find that we have to reduce the timestep slightly, here we use $c\dt/h = 0.4$.

\REV{3}{High-order accurate} Hermite methods were first introduced for hyperbolic systems of equations by Goodrich et al. in~\cite
{goodrich2006hermite}. Since then there have been many enhancements to the original methods described in~\cite
{goodrich2006hermite}. These improvements include but are not limited to order-adaptive implementations~\cite
{chen2012p}, flux-conservative formulation for conservation laws~\cite{kornelus2018flux}, coupling with a discontinuity
sensor to resolve kinks~\cite{alvarez2022hermite} and coupling with discontinuous Galerkin methods to handle complex
boundaries~\cite{beznosov2021hermite,chen2014hybrid}.
In~\cite{appelo2018hermite} the authors developed \textsl{dissipative} (FOT) and \textsl{conservative} (\REV{3}{SOT}) 
Hermite interpolation based schemes for the wave equation. 
The schemes were developed for rectangular geometries where boundary conditions can be imposed
by a simple mirroring strategy. Here these \REV{4}{schemes} are extended to curvilinear geometries.
For complex geometries Hermite methods have been used as efficient building blocks in
hybrid methods. In~\cite{chen2014hybrid} the authors developed a hybrid Hermite-discontinuous Galerkin (DG) scheme for
solving hyperbolic systems and in~\cite{beznosov2021hermite} a Hermite-DG scheme was developed for the wave equation.
The DG approach was used on curvilinear grids near boundaries, but suffered from a smaller time-step restriction than the Hermite method.
One goal of the current work is to develop Hermite schemes for curvilinear domains with boundaries so they
can eventually \REV{4}{be} used on overset grids which consist of overlapping curvilinear grids near boundaries and
one or more background Cartesian grids. In this way efficient high-order accurate Hermite methods with large 
CFL time-steps can be used for complex geometry. \REV{3}{As with all high order accurate methods, to achieve design rates of convergence, Hermite methods require that the solutions are smooth enough. Here we only consider examples that have smooth solutions. Problems with non-smooth solutions discretized with high order accurate Hermite methods were studied in \cite{AppHaglinear} where it was found that (just as for finite difference methods \cite{Banks:2008uq}) they have sublinear rates of convergence but with constants that decrease rapidly with order.}


Compatibility boundary conditions have been used for finite-difference methods for many years (at least since the
early 1980s). For example, in~\cite{henshaw1994fourth,henshaw1994fourthOverlapping,meng2020fourth} the authors use compatibility
conditions for second-order and fourth-order accurate approximations of the incompressible Navier-Stokes equations. For
wave problems, compatibility conditions have been used in~\cite{henshaw2006moving} for compressible Navier-Stokes and
linear elasticity~\cite{appelo2012numerical}, 
as well as high-order schemes for Maxwell's equations~\cite{henshaw2006moving,angel2019high}.
Shu and collaborators have used CBCs in their inverse-Lax-Wendroff approach for
hyperbolic equations and conservation laws~\cite{ILWConservationLaws,ILWConservationLawsSequel,ILWHyperbolicChangingWind,ILWHyperbolic}
as well as for parabolic and advection-diffusion equations~\cite{ILWStabilityParabolic,ILWConvectionDiffusion}.
CBCs are used in the book by Gustafsson 
on high-order difference methods~\cite{GustafssonHighOrderDifferenceMethodsBook2008}.
CBCs have been used by LeVeque and Li with their immersed interface method
to develop accurate approximations at embedded interfaces~\cite{IMElliptic,IMStokesFlow,IMFluidComplexGeometry}.
CBCs have also been used to derive stable and accurate embedded boundary
approximations~\cite{EBWaveOrder2,EBelasticWaveOrder2,EBWaveOrder4}.
CBCs have been incorporated into summation-by-parts schemes by Sj\"ogreen and Petersson for the equations of elasticity~\cite{SjogreenPetersson2012}. 
\REV{4}{Another example where compatibility conditions are used for wave propagation is in the difference potential method by Petropavlovsky et al.~\cite{petropavlovsky2020numerical}.}

In recent work~\cite{hassanieh2021local}, the authors develop~\textsl{local compatibility boundary conditions} (LCBCs) for 
high-order accurate finite difference methods on Cartesian and curvilinear grids. 
The LCBC approach was actually first inspired by the CBC approach described in the present article for use with Hermite methods 
(even though publication of the LCBC method appears first).

A large number of methods for the wave equation  have been proposed in the literature. An incomplete list is, finite
difference methods based on the summation-by-parts framework~\cite{SBP} and upwinding~\cite
{banks2012upwind,Banks2015}, finite element methods which use mass lumping to achieve efficiency~\cite
{Jol-2003}, discontinuous Galerkin methods~\cite{GSSwave,IPDG_Elastic,Upwind2,ChouShuXing2014}, as well as more exotic
methods such as Fourier-Continuation~\cite{FCAD1,FCAD2} and Galerkin differences~\cite{BanksHagstromGD}.  We contend
that Hermite methods have unique properties \REV{5}{(in particular the order independent CFL number in combination with the $\sim 2m^{\rm th}$ order using $\sim (m+1)^d$ degrees of freedoms per element)} that makes their development for use on more general
geometries worthwhile.

The remainder of the paper is organized as follows. 
In Section~\ref{sec:preliminaries} the governing equations are presented together with a discussion
of curvilinear grids and the representation of discrete solutions for the Hermite schemes.
A high-level summary of the Hermite schemes is given in Section~\ref{sec:algorithms}.
Section~\ref{sec:cbc} describes the compatibility boundary condition approach for both Cartesian and
curvilinear grids. An analysis of the solvability and conditioning of the systems of linear equations
resulting from the CBC approximations is given in Section~\ref{sec:AnalysisOfCompat}.
Numerical results are presented in Section~\ref{sec:numericalResults}.
The Appendices contain some proofs, a description of the evolution operators for the Hermite schemes,
as well as some practical implementation details.


\newcommand{\uvBar}{\bar{\uv}}
\newcommand{\uBar}{\bar{u}}
\newcommand{\vvBar}{\bar{\vv}}
\newcommand{\vBar}{\bar{v}}
\section{Preliminaries: governing equations, grids, and discrete approximations} \label{sec:preliminaries}

\subsection{Governing equations}

We consider the initial boundary-value problem for the solution $u=u(\xv,t)$ to 
the scalar wave equation 
\bse
\label{eq:WAVEequationIBVP}
\bat
   & \frac{\p^2 u}{\p t^2} = c^2 \Delta u + f(\xv,t), &&\quad \xv \in \Omega, ~t > 0, \label{eq:WAVEequation} \\
   & u(\xv,0) = U_0(\xv),                             &&\quad \xv \in \Omega ,        \label{eq:waveU0} \\
   & \frac{\p u}{\p t}(\xv) = U_1(\xv) ,              &&\quad \xv \in \Omega, \\
   & \Bc u (\xv,t) = g(\xv,t),                        &&\quad \xv\in\p\Omega.   \label{eq:BC}
\eat
\ese
Here $c>0$ is the wave speed, $\xv=[x_1,x_2]=[x,y]$, and $f(\xv,t)$ and $g(\xv,t)$ are given forcing functions.
The operator $\Bc$ in~\eqref{eq:BC} denotes the boundary condition operator being of Dirichlet or Neumann type.
The initial conditions for $u$ and $\p_t u$ are given in terms of known functions  $U_0(\xv)$ and $U_1(\xv)$.
In this article we focus on problems in $d=2$ space dimensions but note that the method can be extended to three dimensions. 

\subsection{Mappings and parameter space equations}

To discretize~\eqref{eq:WAVEequationIBVP} on a structured curvilinear grid, we assume there exists a smooth and invertible mapping $\Gv$,
\ba
  \xv = \Gv(\rv),
\ea
from the unit square coordinates $\rv=[r_1,r_2]=[r,s]\in[0,1]^2$ to the physical domain coordinates $\xv\in\Omega$.
Using the chain rule the wave equation~\eqref{eq:WAVEequation} can be transformed to parameter space 
\bse
\label{eq:WaveParameterSpace}
\ba
  & \p_t^2 U = c^2 \Delta_\rv U + F(\rv,t), \\
  & c^2 \Delta_\rv U \eqdef  a^{20}(\rv) \frac{\p^2 U}{\p r_1^2} + a^{11}(\rv)\frac{\p^2 U}{\p r_1 \p r_2} + a^{02}(\rv) \frac{\p^2 U}{\p r_2^2} 
         + a^{10}(\rv)\frac{\p U}{\p r_1} + a^{01}(\rv) \frac{\p U }{\p r_2},   
\ea
where $U(\rv,t)=u(\Gv(\rv),t)$, $F(\rv,t)=f(\Gv(\rv),t)$, and the coefficients are 
\ba
   &  a^{20} = c^2 \,  \| \grad_\xv r_1 \|^2 = c^2\, \Big[ (\p_{x_1} r_1)^2 + (\p_{x_2} r_1)^2  \Big], \\
   &  a^{02} = c^2 \,  \| \grad_\xv r_2 \|^2 = c^2\, \Big[ (\p_{x_1} r_2)^2 + (\p_{x_2} r_2)^2  \Big] , \\
   &  a^{11} = c^2 \,  2 (\grad_\xv r_1 \cdot \grad_\xv r_2) = 2(  \p_{x_1} r_1 \, \p_{x_1} r_2 + \p_{x_2} r_1\, \p_{x_2} r_2  ) , \\
   &  a^{10} = c^2 \,  (\p_{x_1}^2 +  \p_{x_2}^2 ) r_1, \\ 
   &  a^{01} = c^2 \,  (\p_{x_1}^2 +  \p_{x_2}^2 ) r_2 .
\ea
\ese
In the remainder of the article we will write $u(\rv,t)$ instead of $U(\rv,t)$. \REV{5}{We note that here we only consider a constant wave speed, but that a spatially variable wave speed can easily be incorporated since we are already discretizing a variable coefficient problem in the reference domain.}

\subsection{Grids, Taylor polynomial representations and Hermite interpolants}

The unit square parameter space will be discretized with a primal (node centered) grid and a dual (cell centered) grid.
Let $\rv_\iv$ denote points on either the primal or dual grid, where $\iv=[i_1,i_2]=[i,j]$ is a multi-index.
The primal grid points are 
\bse
\ba
  &\rv_\iv = [i_1 \dr_1, i_2 \dr_2] , \quad i_k =0,1,2,\ldots,N_k, \\
  & \dr_k = \f{1}{N_k},
\ea
\ese
where $\dr_k$ is the grid spacing in coordinate direction $k$ and $N_k$ is the number of grid cells.
We also use $\dr=\dr_1$ and $\ds=\dr_2$.
The dual grid points are
\ba
   \rv_\iv = [i_1 \dr_1, i_2 \dr_2] , \quad i_k =\half,\f{3}{2},\f{5}{2},\ldots,N_k-\half, 
\ea

The spatial approximation to $u$ near the grid point $\rv_\iv$ is represented as a Taylor polynomial with $(m+1)^{\REV{3}{2}}$ degrees of freedom (DOF),
\bse
\ba
  &  u_\iv(\rv) = \sum_{l_1=0}^m \sum_{l_2=0}^m u_{\iv, l_1,l_2} \, R_i^{l_1} \, S_j^{l_2},  \label{eq:TaylorPoly}\\
  &  R_i \eqdef \f{r - r_i}{\dr} , \quad S_j\eqdef \f{s-s_j}{\ds} ,
\ea
\ese
where the integer $m$ is the degree of the approximation and
where $u_{\iv, l_1,l_2}$ is an approximation to the scaled derivative of $u$,
\ba
   u_{\iv, l_1,l_2} \approx \f{\dr^{l_1}}{l_1 !} \f{\ds^{l_2}}{l_2 !} \p_{r}^{l_1} \p_{s}^{l_2} u(\rv_\iv) .
\ea
We also denote as $u_\iv$ the set of DOFs (or grid function) associated with the Taylor polynomial,
\ba
   u_\iv = \left\{  u_{\iv, l_1,l_2} \right\}_{l_1,l_2=0,1,\ldots,\REV{4}{m}}.
\ea
The Hermite interpolant, centered at $\rv_\iv$, is the polynomial 
that interpolates the solution and it's derivatives at the four neighbouring points $[i_1\pm\half,i_2\pm\half]$
and has the representation, with $(2m+2)^{\REV{3}{2}}$ degrees of freedom, of the form
\ba
  \uBar_\iv(\rv) = \sum_{l_1=0}^{2m+1} \sum_{l_2=0}^{2m+1} \uBar_{\iv, l_1,l_2} \, R_i^{l_1} \, S_j^{l_2} . 
\ea
The over-bar on $\uBar_\iv(\rv)$ will indicate that this representation has $(2m+2)^{\REV{3}{2}}$ DOFs.
See~\ref{sec:HermiteInterpolants} for details on forming the interpolant.

\section{Hermite algorithms}  \label{sec:algorithms}

The basic structure of Hermite scheme (FOT or \REV{3}{SOT}) is given in Algorithm~\ref{alg:hermite}
and illustrated in Figure~\ref{fig:hermiteME2D}.
In the algorithm, $u_\iv^n$ denotes an approximation to the DOFs (solution and derivatives)  at time $t^n = n\dt$, where $\dt$ is the time-step.
The algorithm requires a function $\Ih$ to compute the Hermite interpolant, a function $\Th$ to evolve the solution over a half time-step, and
a function $\Bh$ to assign the boundary conditions.
$P$ denotes the index set of primal points, $D$ the index set of dual points, and $\p P$ the index set of primal boundary points.
Recall that the over-bar on a variable denotes a grid function with $(2m+2)^{\REV{3}{^2}}$ DOFs, while no over-bar is a grid function with $(m+1)^{\REV{3}{^2}}$ DOFs. 
Note that for the FOT scheme, degrees of freedom will be stored for both the solution $u$ and it's time derivative $v=\p_t u$, but
these and other details are left out to simplify the presentation.

{
\newcommand{\algFontSize}{\small}
\begin{algorithm}
\algFontSize 
\caption{Hermite time-stepping algorithm.}
\begin{algorithmic}[1]
  \Function{Hermite}{}  
    \State Compute $\dt$ and number of time-steps $N_t$.
    \State Assign initial conditions.
    \For{$n=1,2,\ldots,N_t$} \Comment Begin time-stepping loop
      \State $t^n = (n-1) \dt$                                                             \Comment Current time.
      \State $\uBar_\jv^n  = \Ih( u_\iv^n )$, \quad $\iv\in P$, $\jv\in D$                 \Comment Interpolate to dual grid.
      \State $u_\jv^{n+\half}  = \Th( \uBar_\jv^n )$, \quad $\jv\in D$                     \Comment Evolve on dual grid to $t^n+\dt/2$.
      \State $\uBar_\kv^{n+\half} = \Bh( u_\jv^{n+\half} )$, \quad $\kv\in\p P$.           \Comment Apply BCs to primal at $t+\dt/2$.
      \State $\uBar_\iv^{n+\half}  = \Ih( u_\jv^{n+\half} )$, \quad $\iv\in P$, $\jv\in D$ \Comment Interpolate to primal. 
      \State $u_\iv^{n+1}  = \Th( \uBar_\iv^{n+\half} )$, \quad $\iv\in P$                 \Comment Evolve on primal to $t^n+\dt$.

    \EndFor    \Comment End time-stepping loop
 \EndFunction

\end{algorithmic} 
\label{alg:hermite}
\end{algorithm}

}

{
\newcommand{\dcol}{red}
\newcommand{\pcol}{blue}
\newcommand{\tikzcircle}[2][red,fill=red]{\tikz[baseline=-0.5ex]\draw[#1,radius=#2] (0,0) circle ;}%
%

\begin{figure}[hbt!]
  \centering
  \begin{tikzpicture}
    \useasboundingbox (0,-1) rectangle (14.25,7.5);

    \begin{scope}[xshift=0cm,yshift=4.cm]
        \draw[black] (1.5,-.2) node[anchor=north] {\smallss 1. Start $t^n$};

       \foreach \y in {0,1,2,3}
         \draw[thick,black!50] (0,\y) -- (3,\y);

       \foreach \x in {0,1,2,3}
         \draw[thick,black!50] (\x,0) -- (\x,3);

        \foreach \y in {0,1,2,3}
          \foreach \x in {0,1,2,3}
           \draw (\x,\y) node {\tikzcircle[blue,very thick, fill=blue!20]{3pt}};

    \end{scope}

    \begin{scope}[xshift=3.75cm,yshift=4.0cm]
        \draw[black] (1.5,-.2) node[anchor=north] {\smallss 2. Interpolate to dual $t^n$};

       \foreach \y in {0,1,2,3}
         \draw[thick,black!50] (0,\y) -- (3,\y);

       \foreach \x in {0,1,2,3}
         \draw[thick,black!50] (\x,0) -- (\x,3);

        \foreach \y in {.5,1.5,2.5}
          \foreach \x in {.5,1.5,2.5}
          {
           \draw (\x,\y) node {\tikzcircle[blue,very thick, fill=blue]{3pt}};
           \draw[<-,thick,blue!50] (\x-.1,\y-.1) -- (\x-.5,\y-.5);
           \draw[<-,thick,blue!50] (\x+.1,\y-.1) -- (\x+.5,\y-.5);
           \draw[<-,thick,blue!50] (\x-.1,\y+.1) -- (\x-.5,\y+.5);
           \draw[<-,thick,blue!50] (\x+.1,\y+.1) -- (\x+.5,\y+.5);
          }

    \end{scope}

    \begin{scope}[xshift=7.5cm,yshift=4.0cm]
        \draw[black] (1.5,-.2) node[anchor=north] {\smallss 3. Evolve dual $t^{n+\half}$};

       \foreach \y in {0,1,2,3}
         \draw[thick,black!50] (0,\y) -- (3,\y);

       \foreach \x in {0,1,2,3}
         \draw[thick,black!50] (\x,0) -- (\x,3);

        \foreach \y in {.5,1.5,2.5}
          \foreach \x in {.5,1.5,2.5}
           \draw (\x,\y) node {\tikzcircle[blue,very thick, fill=blue!20]{3pt}};

    \end{scope}

    \begin{scope}[xshift=11.25cm,yshift=4.0cm]
        \draw[black] (1.5,-.2) node[anchor=north] {\smallss 4. BCs $t^{n+\half}$};

       \foreach \y in {0,1,2,3}
         \draw[thick,black!50] (0,\y) -- (3,\y);

       \foreach \x in {0,1,2,3}
         \draw[thick,black!50] (\x,0) -- (\x,3);

        \foreach \y in {.5,1.5,2.5}
          \foreach \x in {.5,1.5,2.5}
           \draw (\x,\y) node {\tikzcircle[blue,very thick, fill=blue!20]{3pt}};

       \foreach \y in {0,1,2,3}
          \foreach \x in {0,3}
           \draw (\x,\y) node {\tikzcircle[blue,very thick, fill=blue]{3pt}};

      \foreach \y in {0,3}
          \foreach \x in {1,2}
           \draw (\x,\y) node {\tikzcircle[blue,very thick, fill=blue]{3pt}};      

       \foreach \y in {.5,1.5,2.5}
          \foreach \x in {.5}
          {
           \draw[->,thick,blue!50] (\x-.1,\y-.1) -- (\x-.4,\y-.4);
           \draw[->,thick,blue!50] (\x-.1,\y+.1) -- (\x-.4,\y+.4);
          }

      \foreach \y in {.5,1.5,2.5}
          \foreach \x in {2.5}
          {
           \draw[->,thick,blue!50] (\x+.1,\y-.1) -- (\x+.4,\y-.4);
           \draw[->,thick,blue!50] (\x+.1,\y+.1) -- (\x+.4,\y+.4);
          }

       \foreach \y in {.5}
          \foreach \x in {.5,1.5,2.5}
          {
           \draw[->,thick,blue!50] (\x-.1,\y-.1) -- (\x-.4,\y-.4);
           \draw[->,thick,blue!50] (\x+.1,\y-.1) -- (\x+.4,\y-.4);
          }

       \foreach \y in {2.5}
          \foreach \x in {.5,1.5,2.5}
          {
           \draw[->,thick,blue!50] (\x-.1,\y+.1) -- (\x-.4,\y+.4);
           \draw[->,thick,blue!50] (\x+.1,\y+.1) -- (\x+.4,\y+.4);
          }

    \end{scope}

    \begin{scope}[xshift=0cm,yshift=-.2cm]
        \draw[black] (1.5,-.2) node[anchor=north] {\smallss 5. Interp. to primal $t^{n+\half}$};

       \foreach \y in {0,1,2,3}
         \draw[thick,black!50] (0,\y) -- (3,\y);

       \foreach \x in {0,1,2,3}
         \draw[thick,black!50] (\x,0) -- (\x,3);


         \foreach \y in {0,1,2,3}
          \foreach \x in {0,1,2,3}
           \draw (\x,\y) node {\tikzcircle[blue,very thick, fill=blue]{3pt}};

       \foreach \y in {1,2}
          \foreach \x in {1,2}
          {
           \draw[<-,thick,blue!50] (\x-.1,\y-.1) -- (\x-.5,\y-.5);
           \draw[<-,thick,blue!50] (\x+.1,\y-.1) -- (\x+.5,\y-.5);
           \draw[<-,thick,blue!50] (\x-.1,\y+.1) -- (\x-.5,\y+.5);
           \draw[<-,thick,blue!50] (\x+.1,\y+.1) -- (\x+.5,\y+.5);           
          }

    \end{scope}

    \begin{scope}[xshift=3.75cm,yshift=-.2cm]
        \draw[black] (1.5,-.2) node[anchor=north] {\smallss 6. Evolve primal $t^{n+1}$};

       \foreach \y in {0,1,2,3}
         \draw[thick,black!50] (0,\y) -- (3,\y);

       \foreach \x in {0,1,2,3}
         \draw[thick,black!50] (\x,0) -- (\x,3);

        \foreach \y in {0,1,2,3}
          \foreach \x in {0,1,2,3}
           \draw (\x,\y) node {\tikzcircle[blue,very thick, fill=blue!20]{3pt}};

    \end{scope}

    \begin{scope}[xshift=7.5cm,yshift=0.0cm]

       \draw (1,2) node {\tikzcircle[blue,very thick, fill=blue!20]{3pt}};
       \draw[blue] (1,2) node[anchor=west,xshift=5pt] {\smallss $(m+1)^2$ DOF};
       \draw (1,1) node {\tikzcircle[blue,very thick, fill=blue]{3pt}};
       \draw[blue] (1,1) node[anchor=west,xshift=5pt] {\smallss $(2m+2)^2$ DOF};

    \end{scope}

  \end{tikzpicture}
  \caption{Stages of the Hermite FOT and SOT schemes. Open circles have $(m+1)^2$ degrees of freedom while solid circles have $(2m+2)^2$ degrees of freedom.
           }
  \label{fig:hermiteME2D}
\end{figure}
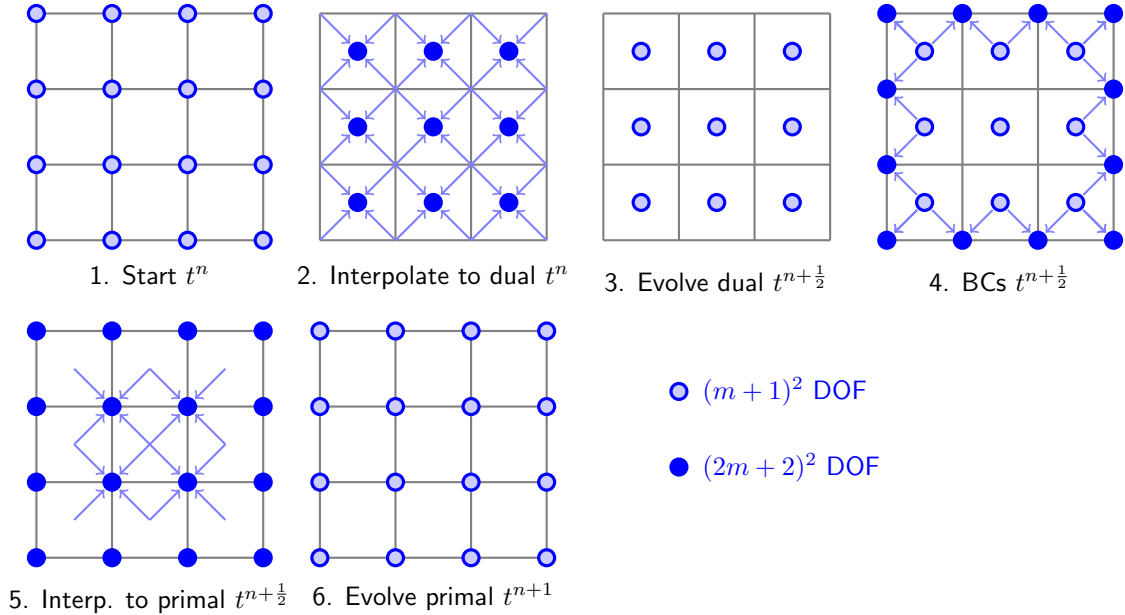 
}

The FOT and \REV{3}{SOT} schemes used here are the extensions, to curvilinear grids, of the \textit{dissipative} 
and \textit{conservative} schemes from~\cite{appelo2018hermite}.
More details of these curvilinear grid schemes are given in~\ref{sec:evolutionAlgorithms} and~\ref{sec:practicalities}.

\REV{3}{
The SOT scheme has less natural dissipation than the FOT scheme. 
When solving problems on curvilinear grids it may be necessary to add some additional dissipation to the SOT scheme so as to retain stability without needing to reduce the time-step. Algorithm~\ref{alg:hermiteSmoothing} shows one way to add dissipation by using the dissipative nature of the Hermite interpolant.
In this algorithm the solution on the primal grid is interpolated to the dual grid, boundary conditions are applied, and then the solution is interpolated back to the
primal. This process can be repeated $N_s$ times.
This smoothing stage can be applied after the second evolution stage in Algorithm~\ref{alg:hermite} (i.e. after line 10).
Smoothing the solution will not change the order of accuracy.
{
\newcommand{\algFontSize}{\small}
\begin{algorithm}
\algFontSize 
\caption{Hermite smoothing for the SOT scheme.}
\begin{algorithmic}[1]
  \Function{HermiteSmooth}{$u_\iv$}  
    \For{$k=1,2,\ldots,N_s$} \Comment Apply $N_s$ smoothing steps
      \State $\uBar_\jv = \Ih( u_\iv )$, \quad $\iv\in P$, $\jv\in D$      \Comment Interpolate to dual grid.
      \State $\uBar_\kv = \Bh( u_\jv )$, \quad $\kv\in\p P$.               \Comment Apply BCs to primal.
      \State $\uBar_\iv = \Ih( u_\jv )$, \quad $\iv\in P$, $\jv\in D$      \Comment Interpolate to primal. 

    \EndFor    \Comment End smoothing steps
 \EndFunction

\end{algorithmic} 
\label{alg:hermiteSmoothing}
\end{algorithm}

}
  
}


\section{Compatibility boundary conditions} \label{sec:cbc}

Compatibility boundary conditions (CBCs) for the wave equation are derived by first taking even time-derivatives of the boundary conditions.
The governing equation~\eqref{eq:WAVEequation} is then used to replace even time-derivatives of $u$ with spatial derivatives of $u$.
The resulting conditions can be used as numerical boundary conditions. These centered conditions are generally more stable and accurate than
using one sided approximations~\cite{hassanieh2021local}.
To illustrate the process for deriving CBCs, consider a Dirichlet boundary condition (Neumann BCs are similar)
\ba
   u(\xv,t)=g(\xv,t), \qquad \xv\in\p\Omega.  \label{eq:Dirichlet}
\ea
Taking two time-derivatives of~\eqref{eq:Dirichlet}
\ba
   \p_t^2 u(\xv,t) = \p_t^2 g(\xv,t),
\ea
and using~\eqref{eq:WAVEequation} gives the first CBC for Dirichlet boundary conditions,
\ba
    L u(\xv,t) + f(\xv,t) = \p_t^2 g(\xv,t), \qquad \xv\in\p\Omega, \label{eq:DirichletCBC1}
\ea
where $L$ is the wave operator,
\ba
   L \eqdef c^2 \Delta.
\ea
The next CBC can be derived by taking two time-derivatives of~\eqref{eq:DirichletCBC1}.
The process can be repeated to derive any number of CBCs.
The CBCs (such as~\eqref{eq:DirichletCBC1}), together with tangential derivatives of the CBCs, 
are used to constrain the Taylor polynomial representation of the solution on the boundary.

\subsection{CBCs on a Cartesian grid}

Figure~\ref{fig:cbc2d} shows a sample grid configuration in two dimensions.
Given values $u_{l_1,l_2,\jv}$ for points $\jv\in D$ on the dual grid (open circles in Figure~\ref{fig:cbc2d}), we require
values $\uBar_{l_1,l_2,\kv}$ for points $\kv\in\p P$ on the boundary of the primal grid (solid circles in Figure~\ref{fig:cbc2d}).

{
\newcommand{\dcol}{red}
\newcommand{\pcol}{blue}
\newcommand{\tikzcircle}[2][red,fill=red]{\tikz[baseline=-0.5ex]\draw[#1,radius=#2] (0,0) circle ;}%
%

\begin{figure}[hbt!]
  \centering
  \begin{tikzpicture}
    \useasboundingbox (0,-.25) rectangle (8,4.25);

    \begin{scope}[xshift=0cm,yshift=0.0cm]

       \foreach \y in {0,1,2,3,4}
         \draw[thick,black!50] (0,\y) -- (4,\y);

       \foreach \x in {0,1,2,3,4}
         \draw[thick,black!50] (\x,0) -- (\x,4);

        \foreach \y in {.5,1.5,2.5,3.5}
          \foreach \x in {.5,1.5,2.5,3.5}
           \draw (\x,\y) node {\tikzcircle[blue,very thick, fill=blue!20]{3pt}};

       \foreach \y in {0,1,2,3,4}
          \foreach \x in {0,4}
           \draw (\x,\y) node {\tikzcircle[blue,very thick, fill=blue]{3pt}};

      \foreach \y in {0,4}
          \foreach \x in {1,2,3}
           \draw (\x,\y) node {\tikzcircle[blue,very thick, fill=blue]{3pt}};      

       \foreach \y in {.5,1.5,2.5,3.5}
          \foreach \x in {.5}
          {
           \draw[->,thick,blue!50] (\x-.1,\y-.1) -- (\x-.4,\y-.4);
           \draw[->,thick,blue!50] (\x-.1,\y+.1) -- (\x-.4,\y+.4);
          }

      \foreach \y in {.5,1.5,2.5,3.5}
          \foreach \x in {3.5}
          {
           \draw[->,thick,blue!50] (\x+.1,\y-.1) -- (\x+.4,\y-.4);
           \draw[->,thick,blue!50] (\x+.1,\y+.1) -- (\x+.4,\y+.4);
          }

       \foreach \y in {.5}
          \foreach \x in {.5,1.5,2.5,3.5}
          {
           \draw[->,thick,blue!50] (\x-.1,\y-.1) -- (\x-.4,\y-.4);
           \draw[->,thick,blue!50] (\x+.1,\y-.1) -- (\x+.4,\y-.4);
          }

       \foreach \y in {3.5}
          \foreach \x in {.5,1.5,2.5,3.5}
          {
           \draw[->,thick,blue!50] (\x-.1,\y+.1) -- (\x-.4,\y+.4);
           \draw[->,thick,blue!50] (\x+.1,\y+.1) -- (\x+.4,\y+.4);
          }

          \draw (0,2) node[anchor=east,xshift=-5pt] {$\xv_{i,j}$};

    \end{scope}

    \begin{scope}[xshift=4.cm,yshift=0.5cm]

       \draw (1,2) node {\tikzcircle[blue,very thick, fill=blue!20]{3pt}};
       \draw[blue] (1,2) node[anchor=west,xshift=5pt] {\smallss $(m+1)^2$ DOF};
       \draw (1,1) node {\tikzcircle[blue,very thick, fill=blue]{3pt}};
       \draw[blue] (1,1) node[anchor=west,xshift=5pt] {\smallss $(2m+2)^2$ DOF};

    \end{scope}

  \end{tikzpicture}
  \caption{Compatibility boundary conditions, together with interior data, are used to define the Hermite representation on the boundary.}
  \label{fig:cbc2d}
\end{figure}
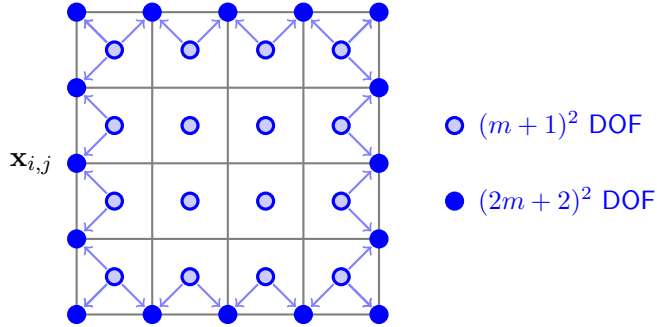 
}












Consider the case of a Cartesian grid for the unit square with grid points $\xv_\iv=(x_i,y_j)$.
Let $\xv_{i,j} = (x_i,y_j)$ be a point on the boundary at $x=0$ (not a corner).
The goal is to define the $(2(m+1))^2$ DOFs in the Taylor polynomial representation of the solution on the boundary, 
\ba 
  & \uBar_{i,j}(\xv) = \sum_{l_1=0}^{2m+1}  \sum_{l_2=0}^{2m+1} \uBar_{l_1,l_2, i,j} \, X_i^{l_1} \, Y_j^{l_2},  \qquad
    X_i \eqdef \f{x - x_i}{\dx} , \quad Y_j\eqdef \f{y-y_j}{\dy} . \label{eq:TaylorPolyCBC}
\ea
As indicated in Figure~\ref{fig:cbc2d}, the polynomial~\eqref{eq:TaylorPolyCBC} is required to match the DOFs from the two nearby interior
points.
The scaled derivatives of the boundary polynomial~\eqref{eq:TaylorPolyCBC} are
\ba
  \f{\dx^\alpha}{\alpha!} \f{\dy^\beta}{\beta!} \, \p_x^\alpha \p_y^\beta \uBar_{i,j}(\xv) = 
    \sum_{l_1=\alpha}^{2m +1}\sum_{l_2=\beta}^{2m +1}  \uBar_{l_1,l_2, i,j} \f{(l_1)(\ldots)(l_1-\alpha+1)}{\alpha ! }\, \f{(l_2)(\ldots)(l_2-\beta+1)}{\beta !}
       \, X^{l_1-\alpha} \, Y^{l_2-\beta}    ,
\ea
and these, when evaluated at $(x_i+\dx/2,y_j\pm\dy/2)$, 
are matched to the given interior dual (scaled) derivatives,
\ba
    u_{\alpha,\beta,i_1+\half,i_2\pm\half} ,    \quad \alpha,\beta=0,1,\ldots,m .
\ea
to give the $2(m+1)^2$ \textsl{interpolation} conditions,
\shadedBoxWithShadow{equation}{orange}{
   \sum_{l_1=\alpha}^{2m +1}\sum_{l_2=\beta}^{2m +1} { l_1 \choose \alpha} \, {l_2 \choose \beta} \,
                 \left[ \half \right]^{l_1-\alpha} \, \left[ \pm\half \right]^{l_2-\beta} \,  \uBar_{l_1,l_2,\iv}   = u_{\alpha,\beta,i_1+\half,i_2\pm\half}, 
                 \quad \alpha,\beta=0,1,\ldots,m .  \label{eq:BoundaryMatchingConditions}
}
A further $2(m+1)^2$ conditions are needed to uniquely determine $\uBar_{l_1,l_2,\iv}$ and these are obtained using CBCs.

\bni
\textbf{Dirichlet boundary conditions.} Consider a Dirichlet boundary condition at $x=0$ 
\ba
   u(0,y,t)=g(y,t). 
\ea
The CBCs are
\ba
   \p_y^\alpha L ^q \, \uBar_{\iv}(0,y_j) = \p_t^{2q} \p_y^\alpha \, g(y_j,t) , \qquad q=0,1,\ldots,m , \quad \alpha=0,1,\ldots,2m+1 .
\ea
Using the binomial expansion gives
\ba
   \p_y^\alpha L^q = c^{2q} \, \p_y^\alpha \, (\p_x^2 + \p_y^2)^q = c^{2q} \, \sum_{k=0}^q {q \choose k} \, \p_x^{2(q-k)} \, \p_y^{2k+\alpha} .
\ea
This leads to the $2(m+1)^2$ CBC conditions
\shadedBoxWithShadow{equation}{orange}{
   c^{2q} \sum_{k=0}^q {q \choose k} \f{(2(q-k))!}{\dx^{2(q-k)}} \, \f{(2k+\alpha)!}{\dy^{2k+\alpha}} \uBar_{2(q-k),2k+\alpha,\iv}
       = \p_t^{2q} \p_y^\alpha g(y_j,t) ,  \label{eq:DirichletHermiteCBCs}
}
for $q=0,1,\ldots,m$, and $\alpha=0,1,\ldots,2m+1$.
Note that in~\eqref{eq:DirichletHermiteCBCs}, only terms with $2k+\alpha\le 2m+1$ should be kept.
Equations~\eqref{eq:BoundaryMatchingConditions} together with~\eqref{eq:DirichletHermiteCBCs} define a linear system of equations
for the unknowns $\uBar_{l_1,l_2,\iv}$ on a Dirichlet boundary.

\bni
\textbf{Neumann boundary conditions.}
Now consider a Neumann boundary condition at $x=0$ ,
\ba
   \p_x u(0,y,t) = g(y,t), 
\ea
The CBCs are
\ba
   \p_y^\alpha \p_x L ^q p(0,y_j) = \p_t^{2q} \p_y^\alpha g(y_j,t) , \qquad q=0,1,\ldots,m , \quad \alpha=0,1,\ldots,2m+1 .
\ea
This leads to the CBC conditions
\shadedBoxWithShadow{equation}{orange}{
   c^{2q} \sum_{k=0}^q {q \choose k} \f{(2(q-k)+1)!}{\dx^{2(q-k)+1}} \, \f{(2k+\alpha)!}{\dy^{2k+\alpha}} \uBar_{2(q-k)+1,2k+\alpha,\iv}
       = \p_t^{2q} \p_y^\alpha g(y_j,t) ,  \label{eq:NeumannHermiteCBCs}
}
for $q=0,1,\ldots,m$, and $\alpha=0,1,\ldots,2m+1$.
Equations~\eqref{eq:BoundaryMatchingConditions} together with~\eqref{eq:NeumannHermiteCBCs} define a linear system of equations
for the unknowns $\uBar_{l_1,l_2,\iv}$ on a Neumann boundary.

\subsection{CBCs for corners on a Cartesian grid} \label{sec:cornerCartesianCBC}

Now consider assigning the solution at the corner point, such as the solution at the lower left point, $\xv_\iv=[0,0]$, in Figure~\ref{fig:cbc2d}.
At this corner there is one interior neighbour on the dual grid at index $\iv=[\half,\half]$. 
Following the discussion for a point on the interior of a side, the conditions to match the coefficients in the Hermite polynomial 
to coefficients of the polynomial at the interior point are the $(m+1)^2$ conditions,
\shadedBoxWithShadow{equation}{orange}{
   \sum_{l_1=\alpha}^{2m +1}\sum_{l_2=\beta}^{2m +1} { l_1 \choose \alpha} \, {l_2 \choose \beta} \,
                 \left[ \half \right]^{l_1-\alpha} \, \left[ \half \right]^{l_2-\beta} \,  \uBar_{l_1,l_2,\iv}   = u_{\alpha,\beta,i_1+\half,i_2+\half}, 
                 \quad \alpha,\beta=0,1,\ldots,m .
                 \label{eq:cornerMatching}
}
A further $3(m+1)^2$ conditions are needed to uniquely determine $\uBar_{l_1,l_2,\iv}$ and these are obtained using CBCs.
In subsequent sub-sections we discuss Dirichlet-Dirichlet (D-D), Neumann-Neumann (N-N), and Dirichlet-Neumann corners (D-N).
This discussion closely follows that given in~\cite{hassanieh2021local} but is adjusted to the setting of Hermite methods.

\subsubsection{Dirichlet-Dirichlet CBC corner}

\mni
Consider a Dirichlet-Dirichlet corner at $\xv=\zerov$,
\bse
\ba
  & u(x,0,t)=g_1(x,t), \\ 
  & u(0,y,t)=g_2(y,t).
\ea
\ese
The CBCs are
\bse
\label{eq:cbcCornerDirichlet}
\ba
   \p_x^\alpha L ^q \uBar_{\iv}(\xv_{\iv}) =  \p_x^\alpha \p_t^{2q} g_1(0,t) ,  \quad q=0,1,\ldots,m,  \\
   \p_y^\alpha L ^q \uBar_{\iv}(\xv_{\iv}) =  \p_y^\alpha \p_t^{2q} g_2(0,t) ,  \quad q=0,1,\ldots,m , 
\ea
\ese
for $\alpha\in \Mc_q$, where $\Mc_q$ is the set of integers,
\ba
  \Mc_q \eqdef \{ 0,1,\ldots,2m +1 \} \, - \, \{ 0,2,4,\dots, 2(q-1) \} .
\ea
The set $\Mc_q$ consists of the integers from 0 to $2m+1$ minus the even integers from $0$ to $2(q-1)$.
In addition, the conditions in~\eqref{eq:cbcCornerDirichlet} should be averaged when $\alpha=2q$.
\REV{1}{The set $\Mc_q$ was derived in~\cite{hassanieh2021local}.
At a corner, the union of the compatibility conditions from adjacent faces are not all independent 
and the set $\Mc_q$ was determined by choosing a independent set of conditions by examining the CBCs in the case of a Cartesian grid.
}
As an example, when $q=1$, $\Mc_1$ is missing $0$, 
\ba
   \Mc_1 = \{ 1,2,3,4,5,\ldots,2m +1 \}, 
\ea
and we average the conditions when $\alpha=2$,
\bse
\bat
  & \p_x L \uBar_{\iv}(\xv_{\iv}) = \p_t^{2q} \p_x \p_t^2 g_1(0,t) ,                                            \qquad&& (\alpha=1), \\ 
  & \p_y L \uBar_{\iv}(\xv_{\iv}) = \p_t^{2q} \p_y \p_t^2 g_2(0,t) ,                                            \qquad&& (\alpha=1), \\ 
  &\half( \p_x^2 L \uBar_{\iv}(\xv_{\iv}) + \p_y^2 L \uBar_{\iv}(\xv_{\iv}) ) = \half ( \p_x^2 \p_t^2 g_1(0,t) + \p_y^2  \p_t^2 g_2(0,t) ), \qquad&& (\alpha=2),  \\
  & \p_x^\alpha L \uBar_{\iv}(\xv_{\iv}) = \p_t^{2q} \p_x^\alpha \p_t^2 g_1(0,t) ,                              \qquad&& \alpha=3,4,5,\ldots,2m+1 \\
  & \p_y^\alpha L \uBar_{\iv}(\xv_{\iv}) = \p_t^{2q} \p_y^\alpha \p_t^2 g_2(0,t) ,                              \qquad&& \alpha=3,4,5,\ldots,2m+1  
\eat
\ese
\REV{1}{
When $q = 2$, $\Mc_2$ is missing $\{0,2\}$,
\ba
	\Mc_2 = \{1,3,4,5,\ldots,2m +1 \},
\ea
and we average the conditions when $\alpha = 4$, 

\bse
\bat
  &\p_x^\alpha L^2 \uBar_{\iv}(\xv_{\iv}) = \p_t^{2q} \p_x^\alpha \p_t^4 g_1(0,t) ,                                            \qquad&& \alpha=1,3, \\ 
  & \p_y^\alpha L^2 \uBar_{\iv}(\xv_{\iv}) = \p_t^{2q} \p_y^\alpha \p_t^4 g_2(0,t)  ,                                            \qquad&& \alpha=1,3, \\ 
  &\half( \p_x^4 L^2 \uBar_{\iv}(\xv_{\iv}) + \p_y^4 L^2 \uBar_{\iv}(\xv_{\iv}) ) = \half ( \p_x^4 \p_t^4 g_1(0,t) + \p_y^4  \p_t^4 g_2(0,t) ), \qquad&& (\alpha=4),  \\
  & \p_x^\alpha L^2 \uBar_{\iv}(\xv_{\iv}) = \p_t^{2q} \p_x^\alpha \p_t^4 g_1(0,t) ,                              \qquad&& \alpha=5,\ldots,2m+1 \\
  & \p_y^\alpha L ^2 \uBar_{\iv}(\xv_{\iv}) = \p_t^{2q} \p_y^\alpha \p_t^4 g_2(0,t) ,                              \qquad&& \alpha=5,\ldots,2m+1  
\eat
\ese
}

\mni
In summary, the CBC D-D corner conditions are (or an average of these conditions when $\alpha=2q$ )
\bse
 \label{eq:CBCcornerDD}
\shadedBoxWithShadow{align}{orange}{
   c^{2q} \sum_{j=0}^q {q \choose j} \f{(2(q-j)+\alpha)!}{\dx^{2(q-j)+\alpha}} \, \f{(2j)!}{\dy^{2j}} \uBar_{2(q-j)+\alpha,2j,\iv}
       = \p_t^{2q} \p_x^\alpha g_1(0,t),  \\
   c^{2q} \sum_{j=0}^q {q \choose j} \f{(2(q-j))!}{\dx^{2(q-j)}} \, \f{(2j+\alpha)!}{\dy^{2j+\alpha}} \uBar_{2(q-j),2j+\alpha,\iv}, 
       = \p_t^{2q} \p_y^\alpha g_2(0,t), 
}
\ese
for $q=0,1,2,\ldots,m$ and $\alpha\in\Mc_q$.
Note that only entries with valid indices $l_1$ and $l_2$ for $\uBar_{l_1,l_2,\iv}$ in \eqref{eq:CBCcornerDD} should be kept.
Equations~\eqref{eq:cornerMatching} and~\eqref{eq:CBCcornerDD} define a linear system of equations to determine the Hermite 
coefficients in the D-D orner at $\xv=\zerov$.

\REV{1}{
Here we count the number of equations to ensure we have a square system. The set $\Mc_q$ has $2m+2 - 2q$ elements. From~\eqref{eq:cbcCornerDirichlet} there are $2(2m+2 - 2q) - 1$ equations obtained, the minus one is a result of averaging when $\alpha = 2q$. Thus, the total number of equations is given by 
\ba
\sum \limits_{q = 0}^{m} 2(2m+2 - 2q) - 1 &= 3(m+1)^2,
\ea 
which is the number of equations required. The arguments for obtaining a square system for the remaining cases follow a similar computation. 
}

\subsubsection{Neumann-Neumann CBC corner}

Now consider a Neumann-Neumann corner at $\xv=\zerov$,
\bse
\ba
  & \p_y u(x,0,t)=g_1(x,t), \\ 
  & \p_x u(0,y,t)=g_2(y,t).
\ea
\ese
The CBCs are
\bse
\label{eq:cbcCornerNeumann}
\ba
   \p_x^\alpha \p_y L ^q \uBar_{\iv}(\xv_{\iv}) =  \p_x^\alpha \p_t^{2q} g_1(0,t) ,  \quad q=0,1,\ldots,m,  \\
   \p_y^\alpha \p_x L ^q \uBar_{\iv}(\xv_{\iv}) =  \p_y^\alpha \p_t^{2q} g_2(0,t) ,  \quad q=0,1,\ldots,m .  
\ea
\ese
Following the argument from the previous section, the CBC N-N corner conditions are thus (or an average of these conditions when $\alpha=2q+1$)
\bse
\label{eq:CBCcornerNN}
\shadedBoxWithShadow{align}{orange}{
   c^{2q} \sum_{j=0}^q {q \choose j} \f{(2(q-j)+\alpha)!}{\dx^{2(q-j)+\alpha}} \, \f{(2j+1)!}{\dy^{2j+1}} \uBar_{2(q-j)+\alpha,2j+1,\iv}
       = \p_t^{2q} \p_x^\alpha g_1(0,t) ,\\
   c^{2q} \sum_{j=0}^q {q \choose j} \f{(2(q-j)+1)!}{\dx^{2(q-j)+1}} \, \f{(2j+\alpha)!}{\dy^{2j+\alpha}} \uBar_{2(q-j)+1,2j+\alpha,\iv}
       = \p_t^{2q} \p_y^\alpha g_2(0,t) ,
}
\ese
for $q=0,1,2,\ldots,2m+1$ and $\alpha\in\Nc_q$.
Here $\Nc_q$ is the set
\ba
  \Nc_q \eqdef \{ 0,1,\ldots,2m +1 \} \, - \, \{ 1,3,5,7,\dots, 2 q-1 \} .
\ea
The set $\Nc_q$ consists of the integers from 0 to $m$ minus the odd integers from $1$ to $2 q-1$.
In addition the conditions in~\eqref{eq:CBCcornerNN} should be averaged when $\alpha=2 q +1$.
\REV{1}{As for the case of a D-D corner, 
the set $\Nc_q$ was derived by choosing an independent set of conditions by examining the CBCs in the case of a Cartesian grid.}
Note that only entries with valid indices $l_1$ and $l_2$ for $\uBar_{l_1,l_2,\iv}$ in \eqref{eq:CBCcornerNN} should be kept.
Equations~\eqref{eq:cornerMatching} and~\eqref{eq:CBCcornerNN} define a linear system of equations to determine the Hermite 
coefficients in the N-N corner at $\xv=\zerov$.

\subsubsection{Dirichlet-Neumann CBC corner}

Consider a Dirichlet-Neumann corner at $\xv=\zerov$, with Neumann on the bottom face, and Dirichlet on the left face, 
\bse
\ba
  &  \p_y u(x,0,t)=g_1(x,t), \\ 
  &  u(0,y,t)=g_2(y,t).
\ea
\ese
The CBCs are
\bse
\label{eq:cbcCornerDirichletNeumannI}
\ba
   \p_x^\alpha \p_y L ^q \uBar_{\iv}(\xv_{\iv}) =  \p_x^\alpha \p_t^{2q} g_1(0,t) ,    \quad q=0,1,\ldots,m,  \\
   \p_y^\beta       L ^q \uBar_{\iv}(\xv_{\iv}) =  \p_y^\beta    \p_t^{2q} g_2(0,t) ,  \quad q=0,1,\ldots,m .
\ea
\ese
Following the argument from the previous section, the CBC D-N corner conditions are thus (or an average of these conditions when $\alpha=2q+1$)
\bse
\label{eq:cbcCornerDirichletNeumann}
\shadedBoxWithShadow{align}{orange}{
   c^{2q} \sum_{j=0}^q {q \choose j} \f{(2(q-j)+\alpha)!}{\dx^{2(q-j)+\alpha}} \, \f{(2j+1)!}{\dy^{2j+1}} \uBar_{2(q-j)+\alpha,2j+1,\iv}
       = \p_t^{2q} \p_x^\alpha g_1(0,t) ,\\
   c^{2q} \sum_{j=0}^q {q \choose j} \f{(2(q-j))!}{\dx^{2(q-j)}} \, \f{(2j+\beta)!}{\dy^{2j+\beta}} \uBar_{2(q-j)+1,2j+\beta,\iv}
       = \p_t^{2q} \p_y^\beta g_2(0,t) ,
}
\ese
for $q=0,1,2,\ldots,m$ and $(\alpha,\beta)\in \Mc_q \times \Nc_q$.
Note that $\alpha$, corresponding to the Neumann BC, is in the set $\Mc_q$ associated with the Dirichlet BC on the face, while
$\beta$, corresponding to the Dirichlet BC is associated with $\Nc_q$.
In addition the conditions in~\eqref{eq:cbcCornerDirichletNeumann} should be averaged when $(\alpha,\beta)=(2 q +1,2 q)$.
Note that only entries with valid indices $l_1$ and $l_2$ for $\uBar_{l_1,l_2,\iv}$ in \eqref{eq:cbcCornerDirichletNeumann} should be kept.
Equations~\eqref{eq:cornerMatching} and~\eqref{eq:cbcCornerDirichletNeumann} define a linear system of equations to determine the Hermite 
coefficients in the D-N corner at $\xv=\zerov$.

\subsection{CBCs on a curvilinear grid}

We now consider the imposition of CBCs on a curvilinear grid. 
As for the Cartesian grid case, the coefficients in the Taylor polynomial representation of the solution for a point on the boundary
will be determined from known interior data together with CBCs. The CBCs become algebraically more complicated on a curvilinear grid and rather than writing a
general formula such as~\eqref{eq:DirichletHermiteCBCs}, a recursion is used to form the equations implied by the CBCs.
%
Let 
\ba
   & \uBar_{\iv}(\rv) = \sum_{l_1=0}^{2m +1} \sum_{l_2=0}^{2m +1} \uBar_{l_1,l_2,\iv} \, R_i^{l_1} \, S_j^{l_2}, \qquad
    R_i \eqdef \f{r-r_i}{\dr},  \quad S_j \eqdef \f{s-s_j}{\ds}, 
\ea
denote the Taylor polynomial representation 
for the solution at a point on the boundary $\rv_{\iv}=(r_i,s_j)$. 
Furthermore, let $L \uBar_{\iv}$ have the Taylor polynomial representation
\ba
   & L \uBar_{\iv}(\rv)  = \sum_{l_1=0}^{2m +1} \sum_{l_2=0}^{2m +1} \dBar_{l_1,l_2,\iv} \, R_i^{l_1} \, S_j^{l_2} , \label{eq:taylorPolyLu}
\ea
\REV{1}{where the application of wave operator $L$ to the Hermite representation is described in~\ref{sec:TaylorForLaplacianInCurvilinear}}.

\subsubsection{Dirichlet CBCs on a curvilinear grid}

\newcommand{\alphat}{\alpha}

The CBCs for a Dirichlet boundary condition at the point $\rv_\iv=(r_i,s_j)$ on the boundary at $r=0$, are {
\ba
  \p_s^{\alphat}  L^q \uBar_{\iv}(\rv_{\iv}) = \p_s^{\alphat}  \p_t^{2q} g(s_j,t), \qquad q=0,1,\ldots,m, \quad \alphat=0,1,\ldots,2m+1 . \label{eq:CBCDirichletCurvilinear}
\ea
For $q=0$ this gives the conditions
\ba
   \f{\alphat!}{\ds^{\alphat}} \uBar_{0,\alphat} =  \p_s^{\alphat}  g(s_j,t), \qquad  \alphat=0,1,\ldots, 2m+1 . \label{eq:dirichletCBC0}
\ea
Let $\uvBar\in \Real^{(2(m+1))^2}$ and $\dvBar\in \Real^{(2(m+1))^2}$ denote the vectors with components $\uBar_{l_1,l_2,\iv}$ and $\dBar_{l_1,l_2,\iv}$, for $l_1,l_2=0,1,\ldots,2m+1$ and fixed $\iv$,
\ba
  & \uvBar    = \begin{bmatrix} \uBar_{0,0,\iv} &  \uBar_{1,0,\iv} &  \uBar_{2,0,\iv} & \ldots & \uBar_{2m+1,2m+1,\iv} \end{bmatrix}^T , \\
  & \dvBar    = \begin{bmatrix} \dBar_{0,0,\iv} &  \dBar_{1,0,\iv} &  \dBar_{2,0,\iv} & \ldots & \dBar_{2m+1,2m+1,\iv} \end{bmatrix}^T , 
\ea
where the dependence of $\uvBar$ and $\dvBar$ on $\iv$ has been suppressed.
We then have the following relationship between the coefficients of $\uBar_{\iv}$ and $L \uBar_{\iv}$, 
\ba
   & \dvBar = \Lh \uvBar,
\ea
where $\Lh$ is the matrix implied by~\eqref{eq:taylorPolyLu} 
(\REV{1}{for details see~\ref{sec:applyingTheWaveOperatorInCurvilinearCoordinates}}).
The CBC~\eqref{eq:dirichletCBC0} for $q=0$ can then be expressed as
\ba
    \f{\alphat!}{\ds^{\alphat}} \ev_{[0,\alphat]}^T \uvBar =  \p_s^{\alphat}  g(s_j,t), \qquad  \alphat=0,1,\ldots,2m+1,
\ea
where $\ev_{[0,\alphat]} = \ev_{\alphat(2m+1)}$ denotes the unit vector corresponding to the entry $\uBar_{0,\alphat}$ in the vector $\uvBar$.
In general we have
\ba
   \f{\alphat!}{\ds^{\alphat}}  \ev_{[0,\alphat]}^T \Lh^q \uvBar =  \p_s^{\alphat} \p_t^{2q} g(s_j,t), \qquad q=0,1,\ldots,m, \quad \alphat=0,1,\ldots,2m+1 .
\ea
The CBCs are thus 
\bse
\label{eq:curvilinearDirichletCBC}
\shadedBoxWithShadow{align}{orange}{
  \f{\alphat!}{\ds^{\alphat}}  \Big[\zv_{\alphat}^{q}\Big]^T \uvBar =  \p_s^{\alphat}  g(s_j,t), \qquad q=0,1,\ldots,m, \quad \alphat=0,1,\ldots,2m+1,
}
where the vectors $\zv_{\alphat}^{q}$ satisfy the recursion
\ba
   & \Big[\zv_{\alphat}^{0}\Big]^T  =  \ev_{[0,\alphat]}^T, \\
   & \Big[\zv_{\alphat}^{q}\Big]^T  =  \Big[\zv_{\alphat}^{q-1}\Big]^T \Lh , \qquad q=1,2,\ldots,m  .
\ea
\ese
The equations in~\eqref{eq:curvilinearDirichletCBC} together with the matching conditions~\eqref{eq:BoundaryMatchingConditions} define
a linear system of equations to determine the unknowns $\uBar_{l_1,l_2,\iv}$.

\REV{4}{
Note that the form of~\eqref{eq:curvilinearDirichletCBC} for a Cartesian grid 
is just the same as~\eqref{eq:curvilinearDirichletCBC} but with $L$ being the Cartesian grid version 
and $\p_n$ and $\p_s$ being replaced by one of $\p_x$ or $\p_y$.
}

\subsubsection{Neumann CBCs on a curvilinear grid}

The CBCs for a Neumann BC at the point $\rv_{\iv}$ on the boundary at $r=0$ are
\bse
\ba
  & \p_s^{\alphat} \p_n   L^q \uBar_{\iv}(\rv_\iv) = \p_s^{\alphat}  \p_t^{2q} g(s_j,t), \qquad q=0,1,\ldots,m, \quad \alphat=0,1,\ldots,2m+1,
\ea
where the normal derivative is
\ba
   \p_n = \nv\cdot\grad = n_1 \p_x + n_2 \p_y = (n_1 r_x + n_2 r_y) \p_r + (n_1 s_x + n_2 s_y) \p_s = b_1(r,s) \p_r + b_2(r,s) \p_s. 
\ea
\ese
Let $\p_n \uBar_{\iv}$ have the Taylor polynomial representation
\ba
   & \p_n \uBar_{\iv}(\rv)  = \sum_{l_1=0}^{2m +1} \sum_{l_2=0}^{2m +1} \nBar_{l_1,l_2,\iv} \, R_i^{l_1} \, S_j^{l_2} . \label{eq:taylorPolyNu}
\ea
Following the previous section, let $\uvBar$ denote the vector of coefficients for 
$\uBar_{l_1,l_2,\iv}$.
In terms of this vector of unknowns, the Neumann CBCs can then be expressed as
\bse
\ba
   \f{\alphat!}{\ds^{\alphat}}   \ev_{[0,\alphat]}^T  \Nh \Lh^q \uvBar =  \p_s^{\alphat}  \p_t^{2q} g(s_j,t) ,
\ea
or
\ba  
    & \f{\alphat!}{\ds^{\alphat}}   \wv_{\alphat}^T \Lh^q \uvBar =   \p_s^{\alphat}  \p_t^{2q} g(s_j,t), \\
    & \wv_{\alphat}^T \eqdef \ev_{[0,\alphat]}^T  \Nh .
\ea
\ese
where $\Nh$ is the matrix implied by~\eqref{eq:taylorPolyNu}.
The vectors 
\bse
\label{eq:curvilinearNeumannCBC}
\ba
   \big[ \zv_{\alphat}^{q} \Big]^T \eqdef \wv_{\alphat}^T \Lh^q ,
\ea
can be computed with the recursion
\ba
   & \Big[\zv_{\alphat}^{0}\Big]^T  =  \ev_{[0,\alphat]}^T \Nh , \\
   & \Big[\zv_{\alphat}^{q}\Big]^T  =  \Big[\zv_{\alphat}^{q-1}\Big]^T \Lh , \qquad q=1,2,\ldots,m  ,
\ea
and the CBCs are
\shadedBoxWithShadow{align}{orange}{
    & \f{\alphat!}{\ds^{\alphat}}  \Big[\zv_{\alphat}^{q}\Big]^T  \uvBar =  \p_s^{\alphat}  \p_t^{2q} g(s_j,t) .
}
\ese
The equations in~\eqref{eq:curvilinearNeumannCBC} together with the matching conditions~\eqref{eq:BoundaryMatchingConditions} define
a linear system of equations to determine the unknowns $\uBar_{l_1,l_2,\iv}$.

\subsubsection{Corner CBCs on a curvilinear grid}

The CBCs at a corner follow the discussion in Section~\ref{sec:cornerCartesianCBC} for a Cartesian grid, except that the equations for the CBCs for the Dirichlet or Neumann case 
on a Cartesian grid should be replaced with the corresponding equations for a curvilinear grid.
\REV{1}{For example, at a D-D corner one should use the interpolation conditions~\eqref{eq:cornerMatching},
and replace the Cartesian grid compatibility conditions~\eqref{eq:cbcCornerDirichlet} with the corresponding curvilinear variants such as~\eqref{eq:CBCDirichletCurvilinear}.
}



\section{Analysis of the CBCs}\label{sec:AnalysisOfCompat}

In this section, we analyze the compatibility boundary conditions developed in Section~\ref{sec:cbc}. 
We consider the solvability and conditioning of the resulting linear system of equations.
We also analyze the symmetry properties that arise when using CBCs on a Cartesian grid. 

We write the boundary conditions developed in Section~\ref{sec:cbc} as a linear system 
\ba
  M \uvBar = \bv, 
\ea
where $\uvBar$ and $\bv$ are vectors of degree $4(m+1)^2$ and $M$ is a matrix of size $4(m+1)^2 \times 4(m+1)^2$. 
This system consists of the interpolation conditions~\eqref{eq:BoundaryMatchingConditions} together with the appropriate CBCs.
In order for this system to be solvable we require the matrix $M$ to be non-singular. 
We aim to derive a set of conditions that guarantees solvability. 


\subsection{Solvability of the CBC matrix systems on Cartesian grids}\label{sec:solvabilityCartesian}

The equations appearing in $M$ should be scaled to improve the conditioning of the matrix. A standard approach to improve the conditioning
is to scale each row by the largest entry in absolute value. This will be called row scaling.
Further improvements to the conditioning can be obtained through a process known as 
equilibration\footnote{The Matlab function \texttt{equilibrate} permutes and rescales the matrix 
to have diagonal entries of magnitude one and off-diagonal entries of magnitude at most one.
The HSL Mathematical Software Library~\cite{scott2023hsl}, \url{http://www.hsl.rl.ac.uk}, also has scaling routines that can be used to improve the conditioning.
}.
For Cartesian grids we consider the max-norm condition number of $M$ defined in the usual way as
\ba
  \kappa_\infty(M) \eqdef \| M \|_\infty  \, \| M^{-1} \|_\infty , 
\ea
where the max-norm of a matrix is the maximum row-sum of the magnitudes of the matrix elements.
We consider, without loss of generality, a boundary at $x=0$.
Let $\tcr$ denote the tall-cell ratio,
\ba
   \tcr \eqdef \f{\dx}{\dy}.
\ea

\begin{theorem}[CBC solvability for Cartesian grids.]
\label{theorem:solvabilityCartesian}
  The matrix $M$ resulting from the CBC equations on a Cartesian grid with Dirichlet or Neumann boundary conditions or at a corner
  where Dirichlet or Neumann boundary conditions meet is nonsingular for $m=1,2,3,4$. 
  Furthermore, the max-norm condition number of the row-scaled $M$ only depends on the ratio $\gamma=\dx/\dy$ 
  and thus remains unchanged as the mesh is refined. 
\end{theorem}
The proof of Theorem~\ref{theorem:solvabilityCartesian} for $m=1,2,3,4$ is given in~\ref{sec:solvabilityCartesianProof}.
We conjecture that the result holds for any positive integer $m$. As an example, here is the form of the condition number 
for $m=1$, for Dirichlet and Neumann boundaries,
\bse
\label{eq:conditionCartesianM1}
\bat
   & \kappa_\infty(M) = \max(41,28+3 \tcr^2)  \times \max(\f{121}{16},1+\tcr^2) , \quad && \text{(Dirichlet)}, \\
   & \kappa_\infty(M) = \max \! \left(\frac{125}{16}, 1+{ \tcr }^{2}, \frac{21}{4}+\frac{3 { \tcr }^{2}}{4}\right) \times \max \! \left(\frac{121}{16}, 1+{ \tcr }^{2}\right),
     \quad && \text{(Neumann)}.
\eat
\ese
Note that the condition numbers in~\eqref{eq:conditionCartesianM1} scale in proportion to $\gamma^2$ and thus become large as $\gamma$ get large. 
This is a known phenomena from other discretizations that use CBCs: the grid spacing normal to the boundary should not be \REV{4}{too} large compared to the grid spacings in the tangential directions.

Table~\ref{tab:conditionNumbersCartesian} summarizes the condition numbers for Dirichlet and Neumann boundaries
and corners for $\tcr=1$. Results for row-scaling and equilibration (Matlab) of the matrix are shown. Equilibration roughly reduces the
condition number by a factor of $10$ over row-scaling.

\begin{table}[hbt]
\begin{center}
\begin{tabular}{|c|c|c|c|c|c|c|c|c|c|c|}
\hline
       & \multicolumn{10}{|c|}{CBC Condition Number $\kappa_\infty(M)$ for Cartesian Grids} \\ \hline
  BC   & \multicolumn{2}{|c|}{$m=1$} & \multicolumn{2}{|c|}{$m=2$} & \multicolumn{2}{|c|}{$m=3$}
       & \multicolumn{2}{|c|}{$m=4$} & \multicolumn{2}{|c|}{$m=5$} \\ \hline
       &  rs      &  eq     &  rs     &   eq    &  rs    &  eq     &  rs     &   eq    &  rs    &  eq     \\
  D    &  3.1e2   &  4.8e1  &  5.8e3  &  7.2e2  &  1.3e5 &  1.2e4  &  3.8e6  &  1.7e5  &  9.7e7 &  2.5e6  \\
  N    &  5.9e1   &  1.4e1  &  7.0e2  &  1.8e2  &  1.4e4 &  3.6e3  &  3.3e5  &  4.4e4  &  7.5e6 &  5.9e5  \\
  D-D  &  5.4e2   &  4.0e1  &  8.3e3  &  5.4e2  &  1.8e5 &  9.0e3  &  4.8e6  &  1.4e5  &  1.2e8 &  2.1e6  \\
  N-N  &  5.7e1   &  7.0e1  &  2.7e2  &  2.6e2  &  2.9e3 &  6.7e2  &  3.8e4  &  1.1e4  &  6.3e5 &  1.7e5  \\
  D-N  &  1.6e2   &  1.8e1  &  1.4e3  &  5.0e1  &  2.2e4 &  2.9e3  &  4.2e5  &  1.4e5  &  8.5e6 &  7.8e5  \\ 
  \hline
\end{tabular} 
\end{center}  
\caption{Max-norm condition numbers of the CBC matrices with row-scaling (rs) and equilibration (eq) for $\tcr=\dx/\dy=1$. 
D denotes a Dirichlet BC, N a Neumann BC, D-D a Dirichlet-Dirichlet corner and so on.
  } 
\label{tab:conditionNumbersCartesian}
\end{table}

\subsection{Solvability of CBC matrix systems on curvilinear grids}\label{sec:solvabilityCurvilinear}

To study the solvability of the CBC matrix on a curvilinear grid we freeze coefficients near a point on the boundary and
consider the wave equation with the constant coefficient operator 
\begin{equation}
   L_0 = c^{20} \, \p_{r}^2 + 2 c^{11}  \, \p_r \p_s + c^{02} \, \p_{s}^2  + c^{10} \, \p_r  + c^{01} \, \p_s.\label{eq:constantCoeffOperator}
\end{equation}
The matrix $M$ in the CBC matrix can be formed symbolically for this case and it's determinant can be evaluated.
This leads to the following result.

\begin{theorem}[CBC solvability for Curvilinear grids.]
The CBC matrix $M$ for the constant-coefficient operator $L_0$ in~\eqref{eq:constantCoeffOperator} 
with Dirichlet or Neumann boundary conditions, for $m = 1,2,3,4$, 
is non-singular provided $c^{20} > 0$ and $\dr$ is sufficiently small (left or right face) 
or $c^{02} > 0$ and $\ds$ is sufficiently small (top or bottom face). 
If $c^{10} = 0$ (left or right face) or $c^{01} = 0$ (top or bottom face), 
then the matrix is non-singular for any $\dr$ and $\ds$.
\end{theorem}

\begin{proof}
We will focus on the left boundary at $r=0$, the arguments for the other boundaries are done similarly. 
For Dirichlet or Neumann boundary conditions, the determinant of the matrix $M$,  has the form 
\begin{equation}
\det(M) = K_m G_m(\xi), \quad \xi \eqdef \frac{c^{10}\dr}{2 c^{20}}, \quad m = 1,2,3,4,
\label{eq:LeftSideDet}
\end{equation}
where $K_m$ is a non-zero constant that depends on $\dr$, $\ds$ and $c^{20}$, and $G_m(\xi)$ is a polynomial with $G_m(0) = 1$. 
For Dirichlet boundary conditions the polynomials are given by 
\bse
\ba
G_1(\xi) &=  \left[1 - \frac{\xi}{4}\right]^4,\\
G_2(\xi) &=  \left[ 1-\frac{9}{16} \xi +\frac{3}{32} \xi^{2}-\frac{1}{192} \xi^{3} \right]^6 , \\
G_3(\xi) &=  \left[ 1-\frac{29}{32} \xi +\frac{5}{16} \xi^{2}-\frac{7}{128} \xi^{3}+\frac{1}{192} \xi^{4}-\frac{1}{3840} \xi^{5}+\frac{1}{184320} \xi^{6} \right]^8 ,\\
G_4(\xi) &=  \Big[ 1-\frac{325}{256} \xi +\frac{345}{512} \xi^{2}-\frac{155}{768} \xi^{3}+\frac{235}{6144} \xi^{4}-\frac{99}{20480} \xi^{5}+\frac{245}{589824} \xi^{6}\\
         &\qquad -\frac{199}{8257536} \xi^{7}+\frac{5}{5505024} \xi^{8}-\frac{1}{49545216} \xi^{9}+\frac{1}{4954521600} \xi^{10} \Big]^{10}.
\ea
\ese
For Neumann boundary conditions the corresponding polynomials are
\bse
\ba
    G_1(\xi) &= \left[ 1 - \f{\xi}{2} \right ]^{4}, \\
    G_2(\xi) &= \Big[ 1-\xi +\frac{1}{4} \xi^{2}-\frac{1}{48} \xi^{3} \Big]^6  , \\
    G_3(\xi) &= \Big[ 1-\frac{3}{2} \xi +\frac{3}{4} \xi^{2}-\frac{35}{192} \xi^{3}+\frac{3}{128} \xi^{4}-\frac{1}{640} \xi^{5}+\frac{1}{23040} \xi^{6}
     \Big]^8 , \\
     G_4(\xi) &= \Big[
1-2 \xi +\frac{3}{2} \xi^{2}-\frac{233}{384} \xi^{3}+\frac{29}{192} \xi^{4}-\frac{63}{2560} \xi^{5}+\frac{31}{11520} \xi^{6} \\
       &\qquad -\frac{127}{645120} 
               \xi^{7}+\frac{1}{107520} \xi^{8}-\frac{1}{3870720} \xi^{9}+\frac{1}{309657600} \xi^{10}   \Big]^{10}  . 
                \label{eq:NeumannCBCdet4}
\ea
\ese
For $\xi$ sufficiently small it can be seen that $G_m(\xi)>0$ and thus the matrix $M$ is nonsingular. 
\REV{1}{In particular, Table~\ref{tab:maxXiValues} gives the largest $\xi=\xi_{\rm max}$ such that $G_m(\xi) > 0$ for $0\le \xi \le \xi_{\rm max}$ for the functions above.}
\qed
\end{proof}

\REV{1}{
\begin{table}[hbt]
\begin{center}
\begin{tabular}{|c|c|c|c|c|}
\cline{2-5}
   \multicolumn{1}{c|}{} &  \multicolumn{4}{c|}{ $\xi_{\rm max}$ }   \\
\hline
  BC &  $G_1$  & $G_2$     &  $G_3$ & $G_4$  \\
\hline
Dirichlet &  4.00  & 2.23  & 2.65 & 2.35 \\
Neumann &  2.00  & 1.48  & 2.23 & 1.08 \\
\hline
\end{tabular} 
\end{center}  
\caption{ Largest $\xi=\xi_{\rm max}$ such that $G_m(\xi) > 0$ for $0\le \xi \le \xi_{\rm max}$ for the Dirichlet and Neumann case.
  } 
\label{tab:maxXiValues}
\end{table}
}

To study the actual conditioning of the CBC matrices in practice,
Figure~\ref{fig:cbcMatrixConditionNumbers} graphs $\kappa_{\infty}(M)$ for the CBC matrices on two curvilinear grids as the mesh is refined.
The polynomial mapping and the X mapping \REV{1}{defined in Section~\ref{sec:mappingsAndGrids} } are considered (see Figure~\ref{fig:gridPlots}) for different boundary and corner conditions.
Results are shown for the row scaled matrix and the equilibrated matrix.
The condition numbers are seen to increase with the Hermite degree $m$. The condition numbers 
are fairly constant or sometimes decrease as $\dr$ decreases.
{{}
\newcommand{\figw}{5.25cm}
\newcommand{\figh}{5.25cm}
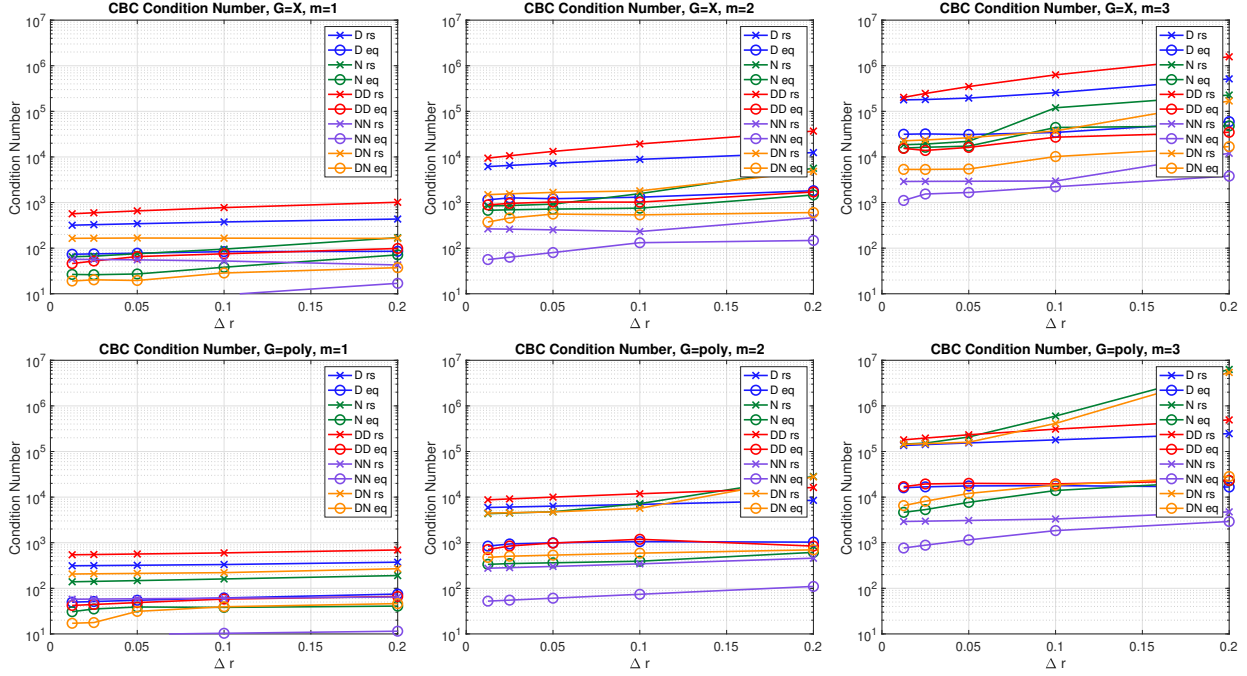
\begin{figure}[htb]
\begin{center}
\begin{tikzpicture}
  \useasboundingbox (0,.5) rectangle (16.5,9);  

  \begin{scope}[yshift=4.5cm]
     \figByWidth{0.0}{0}{fig/hermite2dGXM1ConditionNumber}{\figw}[0][0][0][0];
     \figByWidth{5.5}{0}{fig/hermite2dGXM2ConditionNumber}{\figw}[0][0][0][0];
     \figByWidth{11.}{0}{fig/hermite2dGXM3ConditionNumber}{\figw}[0][0][0][0];
  \end{scope} 
  \begin{scope}[yshift=0cm]
    \figByWidth{0.0}{0}{fig/hermite2dGpolyM1ConditionNumber}{\figw}[0][0][0][0];
    \figByWidth{5.5}{0}{fig/hermite2dGpolyM2ConditionNumber}{\figw}[0][0][0][0];
   \figByWidth{11.}{0}{fig/hermite2dGpolyM3ConditionNumber}{\figw}[0][0][0][0];  
  \end{scope}   
\end{tikzpicture}
\end{center}
\caption{
CBC matrix condition numbers on curvilinear grids using row-scaling (rs) and equilibration (eq) for a Dirichlet BC (D) , Neumann (N) BC, D-D corner, N-N corner, and D-N corner. 
Top: polynomial mapping.
  Bottom: X mapping. Left column: $m=1$. Middle column: $m=2$. Right column: $m=3$.
 }
\label{fig:cbcMatrixConditionNumbers}
\end{figure}
}

\subsection{Symmetry properties of the CBC conditions}\label{sec:symmetrytyOfCompat}

Consider the case of a Cartesian grid with homogeneous Dirichlet or Neumann boundary conditions.
Let us focus on the boundary at $x=0$, the result for other boundaries will be similar.
The CBCs are
\ba
   \Delta^n \, \p_x u(0,y) = [ \p_x^2 + \p_y^2]^n \, \p_x^\alpha u(0,y) = 0 , \quad n=0,1,2,\ldots , 
\ea
where $\alpha=0$ for Dirichlet boundary conditions and $\alpha=1$ for Neumann.
It follows from the binomial expansion that 
 \ba
   \p_x^{2n} \, \p_x^\alpha u(0,y)  = - \sum_{j=1}^n {n \choose j} (\p_x^2)^{n-j} (\p_y^2)^j \, \p_x^\alpha u(0,y). \label{eq:binomialCBC}
\ea
For $n=0$ we have 
\ba
 \p_x^\alpha u(0,y)=0, 
\ea
which implies $\p_y^\beta \p_x^\alpha u(0,y)=0$ for $\beta=0,1,2,\ldots$.
Using this in the right-hand-side of~\eqref{eq:binomialCBC} for $n=1$ gives
\ba
  \p_x^2 \p_x^\alpha(0,y) = 0 , 
\ea
which in turn implies $\p_y^\beta \p_x^2 \p_x^\alpha(0,y)=0$ for $\beta=0,1,2,\ldots$.
This can be used in~\eqref{eq:binomialCBC} for $n=2$ to show $\p_x^4 \p_x^\alpha(0,y)=0$.
Repeating this argument leads to 
\ba
  \p_x^{2n} \p_x^\alpha(0,y) = 0 , \quad n=0,1,2,\ldots
\ea
On a Dirichlet boundary with $u(0,y,t)=0$ it then follows that all even $x$-derivatives of $u$ are zero on the boundary,
\ba
    \p_x^{2n} u(0,y,t) = 0 , \qquad n=0,1,2,\ldots. \label{eq:oddSymmetry}
\ea
This implies that $u$ has odd symmetry in $x$ at the boundary.
On a Neumann boundary with $\p_x u(0,y,t)=0$ it follows that all odd $x$-derivatives of $u$ are zero
\ba
   \p_x^{2n+1} u(0,y,t) = 0 , \qquad n=0,1,2,\ldots ,  \label{eq:evenymmetry}
\ea
and the solution has even symmetry in $x$ at the boundary.
The conditions~\eqref{eq:oddSymmetry} and~\eqref{eq:evenymmetry} are often used as a simple way to set numerical boundary conditions
by odd or even reflection.

The next Theorem shows that the CBC approach leads to a Taylor polynomial representation that has these same symmetry conditions. 
\begin{theorem}[Symmetry of the CBC conditions]
  The Taylor polynomial representation of the solution, resulting from application of the CBC conditions
  on the boundary at $x=0$ of a Cartesian grid,
   has odd symmetry
  for homogeneous Dirichlet boundary conditions and even symmetry for homogeneous Neumann conditions.
\end{theorem}
\begin{proof}
The Taylor polynomial representation on the boundary $x_i=0$ ($i=0$) is 
\ba 
  & \uBar_{i,j}(\xv) = \sum_{l_1=0}^{2m+1}  \sum_{l_2=0}^{2m+1} \uBar_{l_1,l_2, i,j} \, X_i^{l_1} \, Y_j^{l_2}, 
   \qquad X_i \eqdef \f{x}{\dx} , \quad Y_j\eqdef \f{y-y_j}{\dy} . \label{eq:TaylorPolyCBCII}
\ea
For a homogeneous Dirichlet boundary condition the CBC conditions that are imposed are 
\ba
   \p_y^\alpha  \Delta^q u(0,y_j) = 0 , \quad q=0,1,2,\ldots,m, \quad \alpha=0,1,2,\ldots,2m+1.
\ea
Following the argument leading to~\eqref{eq:oddSymmetry} these conditions imply
\ba
   \p_y^\alpha \p_x^{2q} u(0,y_j) = 0 , \quad q=0,1,2,\ldots,m, \quad \alpha=0,1,2,\ldots,2m+1.
\ea
Whence
\ba
   \uBar_{l_1,l_2,i,j} =0 , \quad l_1=0,2,4,\ldots,2m, \quad l_2=0,1,2,\ldots,2m+1.
\ea
Therefore only odd powers of $x$ remain in the Taylor polynomial, which gives the desired result.
For example, at $y=y_j$, the polynomial takes the form
\ba
  \uBar_{i,j}(x,y_j) = \uBar_{1,0,i,j} \, x + \uBar_{3,0,i,j} \, x^3 + \ldots \uBar_{2m+1,0,i,j} \, x^{2m+1} .
\ea
The result for Neumann boundary conditions follows by a similar argument.
\qed
\end{proof}

\newcommand{\kvHat}{\hat{\kv}}
\newcommand{\nTheta}{n_\theta}
\section{Numerical results} \label{sec:numericalResults}

Numerical results are now presented to demonstrate the accuracy and stability of the Hermite schemes on
curvilinear grids with compatibility boundary conditions. Results are shown for both orthogonal and
non-orthogonal grids using Dirichlet and Neumann boundary conditions.

\subsection{Mappings and Grids} \label{sec:mappingsAndGrids}

Plots of the grids used to evaluate the Hermite schemes are shown in Figure~\ref{fig:gridPlots}.
Two of the grids, the rhombus and X mapping, are non-orthogonal. The grids are defined in terms of 
mappings as defined next.
{
\newcommand{\figw}{5.5cm}
\newcommand{\figh}{5.5cm}
\begin{figure}[htb]
\begin{center}
\begin{tikzpicture}
  \useasboundingbox (0,.65) rectangle (16,9.5);  

  \begin{scope}[yshift=4.9cm]
     \figByWidth{0}{0}{fig/GridPlotidentity}{\figw}[0][0][0][0];
     \figByWidth{5.5}{0}{fig/GridPlotpoly}{\figw}[0][0][0][0];
     \figByWidth{11.}{0}{fig/GridPlottanh}{\figw}[0][0][0][0];
  \end{scope} 
  \begin{scope}[yshift=0cm]
    \figByWidth{0}{0}{fig/GridPlotannulus}{\figw}[0][0][0][0];
    \figByWidth{5.5}{0}{fig/GridPlotrhombus}{\figw}[0][0][0][0];
    \figByWidth{11}{0}{fig/GridPlotX}{\figw}[0][0][0][0];
  \end{scope}   
\end{tikzpicture}
\end{center}
\caption{
Plots of the grids used in testing the Hermite schemes.
Top, left to right: identity, polynomial and tanh grids.
Bottom, left to right: annulus, rhombus and X grids.
 }
\label{fig:gridPlots}
\end{figure}
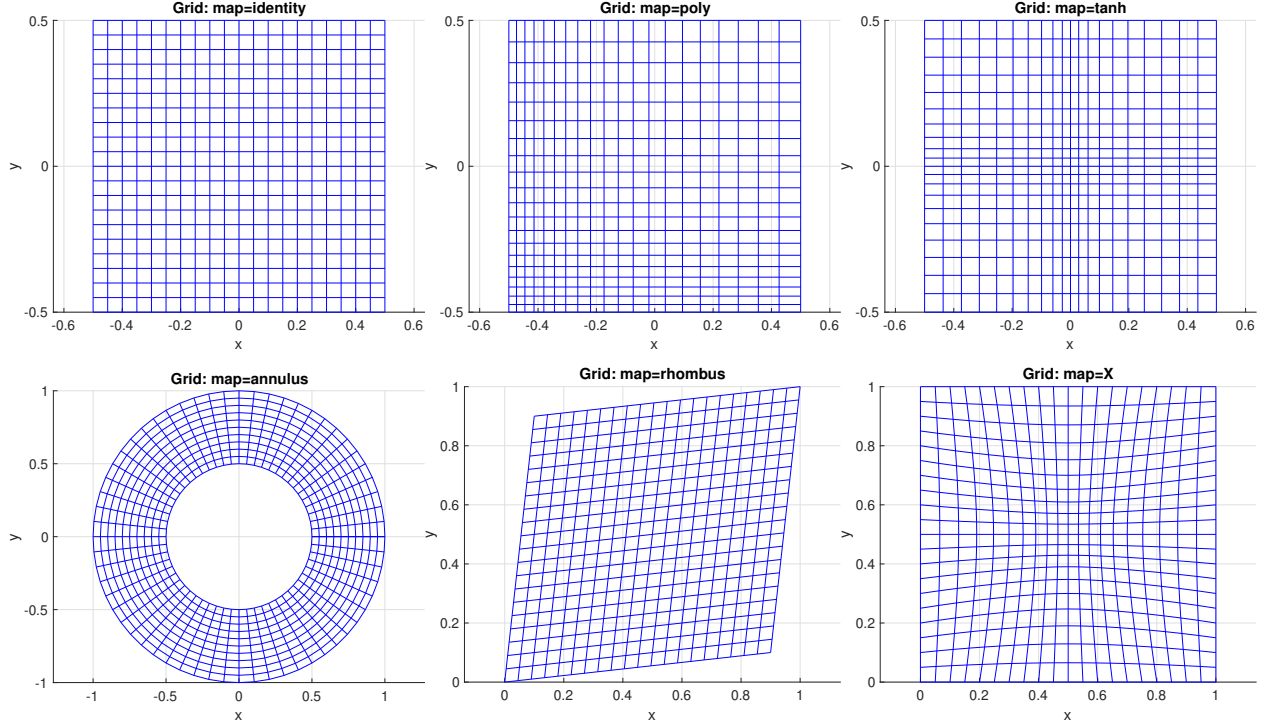
}

\mni
\textbf{Polynomial mapping.}
The \textit{polynomial} mapping can be used to cluster points near a boundary.
In one-dimension it takes the form 
\ba
   x = G(r) = x_a + (x_b-x_a) \Big( \alpha r + (1-\alpha) r^2 \Big) ,
\label{eq:polyMapping}
\ea
and maps $r\in[0,1]$ to $x\in[x_a,x_b]$.
We choose $\alpha=0.5$. This mapping is applied in both the $r_1$ and $r_2$ directions 
to give the polynomial grid in Figure~\ref{fig:gridPlots}.

\mni
\textbf{Hyperbolic tangent mapping.}
The \textit{tanh} mapping can be used to cluster points in the interior of the domain.
In one dimension it is defined by
\bse
\ba
  & x = G(r) =  x_a + (x_b-x_a) \Big( \alpha r + a ( \tanh(\beta (r-r_0)) - \tanh(\beta(-r_0)) ) \Big), \\
  & \alpha = 1 - a \, \big( \tanh(\beta (1-r_0)) - \tanh(\beta (-r_0)) \big) , 
\ea
\ese
where $a$ is an amplitude and $\alpha$ is chosen so $x(1)=x_b$.
We take $r_0=0.5$, $\beta=5$, and $a=-0.15$.

\mni
\textbf{Rhombus mapping.}
The rhombus mapping is a simple non-orthogonal mapping defined by 
\ba
    \xv = \Gv(\rv) 
          = \begin{bmatrix}
                (1-\alpha) r_1 + \alpha r_2 \\
                (1-\beta) r_2 + \beta r_1
            \end{bmatrix} , 
\ea
where we choose $\alpha=0.1$ and $\beta=0.1$

\mni
\textbf{X mapping.}
The X mapping is a non-orthogonal mapping defined by 
\ba
    \xv = \Gv(\rv) 
          = \begin{bmatrix}
                r_1 + \beta r_2(1-r_2) \sin(2 \pi r_1) \\
                r_2 + \beta r_1(1-r_1) \sin(2 \pi r_2) 
            \end{bmatrix} ,
\ea
where we take $\beta=0.2$.


\subsection{Manufactured and exact solutions} \label{sec:manufactured and exact solutions}

{
\newcommand{\figw}{6.5cm}
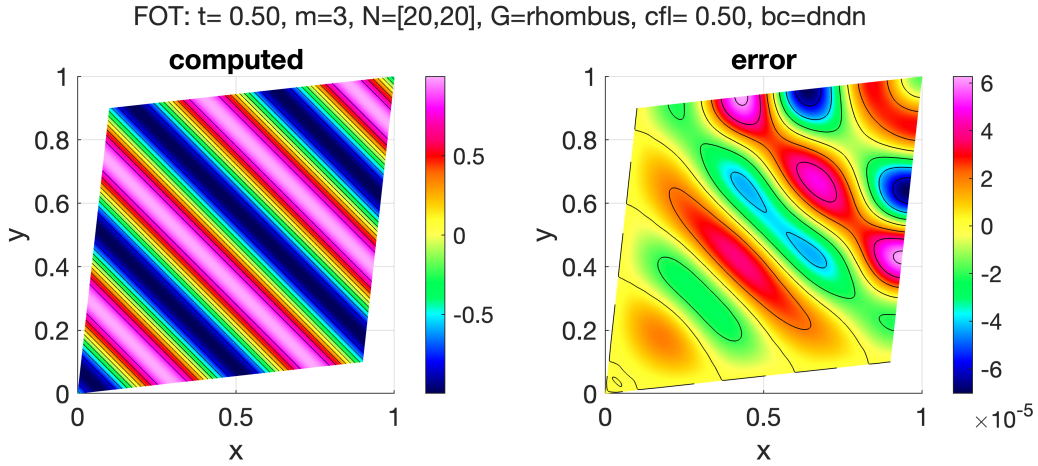
\begin{figure}[htb]
\begin{center}
\begin{tikzpicture}
  \useasboundingbox (0,.85) rectangle (14,6.25);  

  \begin{scope}[yshift=0cm]
    \figByHeight{0}{0}{fig/sineFOT3MapRhombusSolutionAndError}{\figw}[0.05][0][0][0];
    
  \end{scope}   
\end{tikzpicture}
\end{center}
\caption{
Rhombus. Computed solution and error using the FOT scheme with $m=3$ (order $5$) and the sine solution. Boundary conditions are Dirichlet (left, bottom) and
Neumann (right and top).
 }
\label{fig:rhombusSineFOT}
\end{figure}
}

\mni
\textbf{Sine solution.}
The sine solution (shown in Figure~\ref{fig:rhombusSineFOT}) is 
\bse
\label{eq:sineSolution}
\ba
   & u(\xv,t) = \sin( k_x x + k_y y  - \omega t ), \\
   & \omega  = c \, \sqrt{ k_x^2 + k_y^2}.
\ea
\ese
This is an exact solution to the free space problem but requires inhomogeneous boundary conditions.

{
\newcommand{\figw}{6.5cm}
\begin{figure}[htb]
\begin{center}
\begin{tikzpicture}
  \useasboundingbox (0,.85) rectangle (14,6.25);  

  \begin{scope}[yshift=0cm]
    \figByHeight{0}{0}{fig/squareEigMapXSolutionAndError}{\figw}[.05][0][0][0];
  \end{scope}   
\end{tikzpicture}
\end{center}
\caption{
Square eigenfunction on the X grid. Computed solution and error using the SOT scheme with $m=4$ (order 8). 
 }
\label{fig:squareEigGridX}
\end{figure}
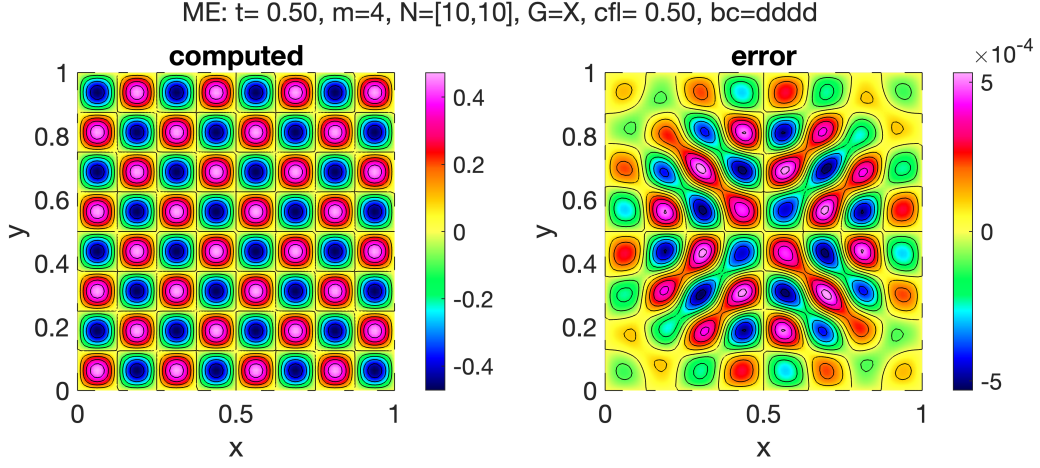
}

\mni
\textbf{Eigenfunction of a square.}
Eigenfunctions of the unit square with Dirichlet boundary conditions (see Figure~\ref{fig:squareEigGridX}) take the form
\bse
\label{eq:sqaureEig}
\ba
    u(\xv,t) = \sin( \pi k_x x ) \sin( \pi k_y y ) \cos( \omega t ) ,
\ea
for integer values of $k_x$ and $k_y$, where
\ba
   \omega = c \, \sqrt{ (\pi k_x)^2 + (\pi k_y)^2 }.
\ea
\ese
Similar expressions can be found for Neumann boundary conditions or a mix of Dirichlet and Neumann boundary conditions.

{
\newcommand{\figw}{6.5cm}
\begin{figure}[htb]
\begin{center}
\begin{tikzpicture}
  \useasboundingbox (0,.85) rectangle (14,6.25);  

  \begin{scope}[yshift=0cm]
    \figByHeight{0}{0}{fig/annulusEigME5n30m1SolutionAndError}{\figw}[0.05][0][0][0];
  \end{scope}   
\end{tikzpicture}
\end{center}
\caption{
Annulus eigenfunction. Computed solution and error using the SOT scheme with $m=5$ (order $10$).
 }
\label{fig:annulusEigME}
\end{figure}
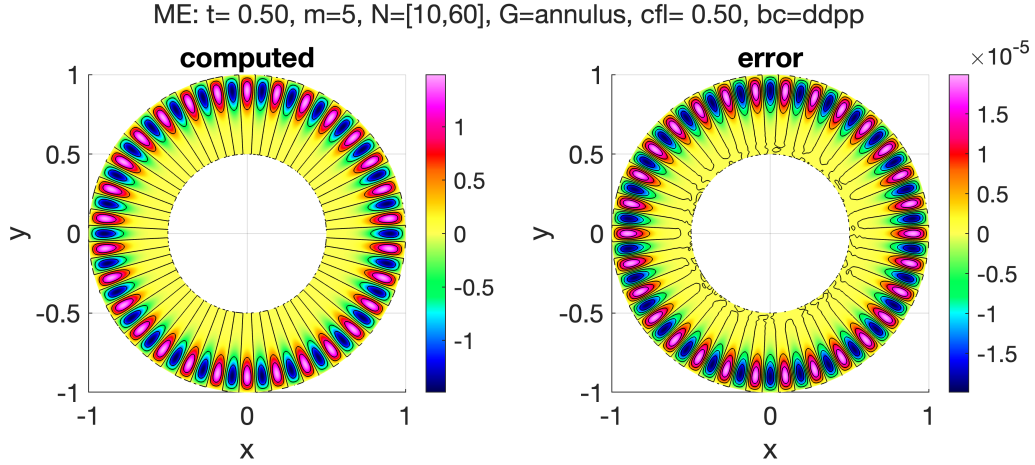
}

\newcommand{\lam}{\lambda_{\nTheta,n_r}}
\textbf{Eigenfunction of an annulus.}
Eigenfunctions of an annulus with an inner radius $r_a=0.5$ and outer radius $r_b=1.0$ and with Dirichlet boundary conditions,
as shown in Figure~\ref{fig:annulusEigME}, are of the form
\bse
\label{eq:annulusEig}
\ba
  &  u(r,\theta,t) = \f{1}{\sqrt{ c_J^2 + c_Y^2}} \Big( c_J J_{\nTheta}(\lam r) + c_Y Y_{\nTheta} (\lam r) \Big) \, \cos( \nTheta \theta ) \cos( c \lam t), \\
  &    c_J = Y_{\nTheta}(\lam r_a), \qquad 
       c_Y = -J_{\nTheta}(\lam r_a),
\ea
\ese
where $J_{\nTheta}$ and $Y_{\nTheta}$ are the Bessel functions of the first kind.
The values of the eigenvalues $\lam$ are 
roots of $d(\lambda)= J_{\nTheta}(\lambda r_a) Y_{\nTheta}(\lambda r_b) - J_{\nTheta}(\lambda r_b) Y_{\nTheta}(\lambda r_a)=0$.

\subsection{Results for the first-order in time (FOT) scheme} \label{sec:fotResults}

{
\newcommand{\figw}{5.5cm}
\newcommand{\figh}{5.5cm}
\begin{figure}[htb]
\begin{center}
\begin{tikzpicture}
  \useasboundingbox (0,.5) rectangle (16.5,12);  

  \begin{scope}[yshift=6cm]
     \figByWidth{0.0}{0}{fig/FOT2DsquareEigBCdMaptConvergence}{\figw}[0][0][0][0];
     \figByWidth{5.5}{0}{fig/FOT2DannulusEigBCdMapaConvergence}{\figw}[0][0][0][0];
     \figByWidth{11.}{0}{fig/FOT2DsineBCdMaptConvergence}{\figw}[0][0][0][0];
  \end{scope} 
  \begin{scope}[yshift=0cm]
    \figByWidth{0}{0}{fig/FOT2DsineBCnMappConvergence}{\figw}[0][0][0][0];
    \figByWidth{5.5}{0}{fig/FOT2DsineBCdndnMaptConvergence}{\figw}[0][0][0][0];
    \figByWidth{11.}{0}{fig/FOT2DsineBCdMapxConvergence}{\figw}[0][0][0][0];
  \end{scope}   
\end{tikzpicture}
\end{center}
\caption{
Grid convergence, Hermite FOT scheme in 2D. Figure titles give run details: S=solution type, B=boundary conditions, G=grid type.
The expected order of accuracy for the FOT scheme is $2m-1$ for a degree $m$ Hermite approximation.
 }
\label{fig:hermite2DFOT}
\end{figure}
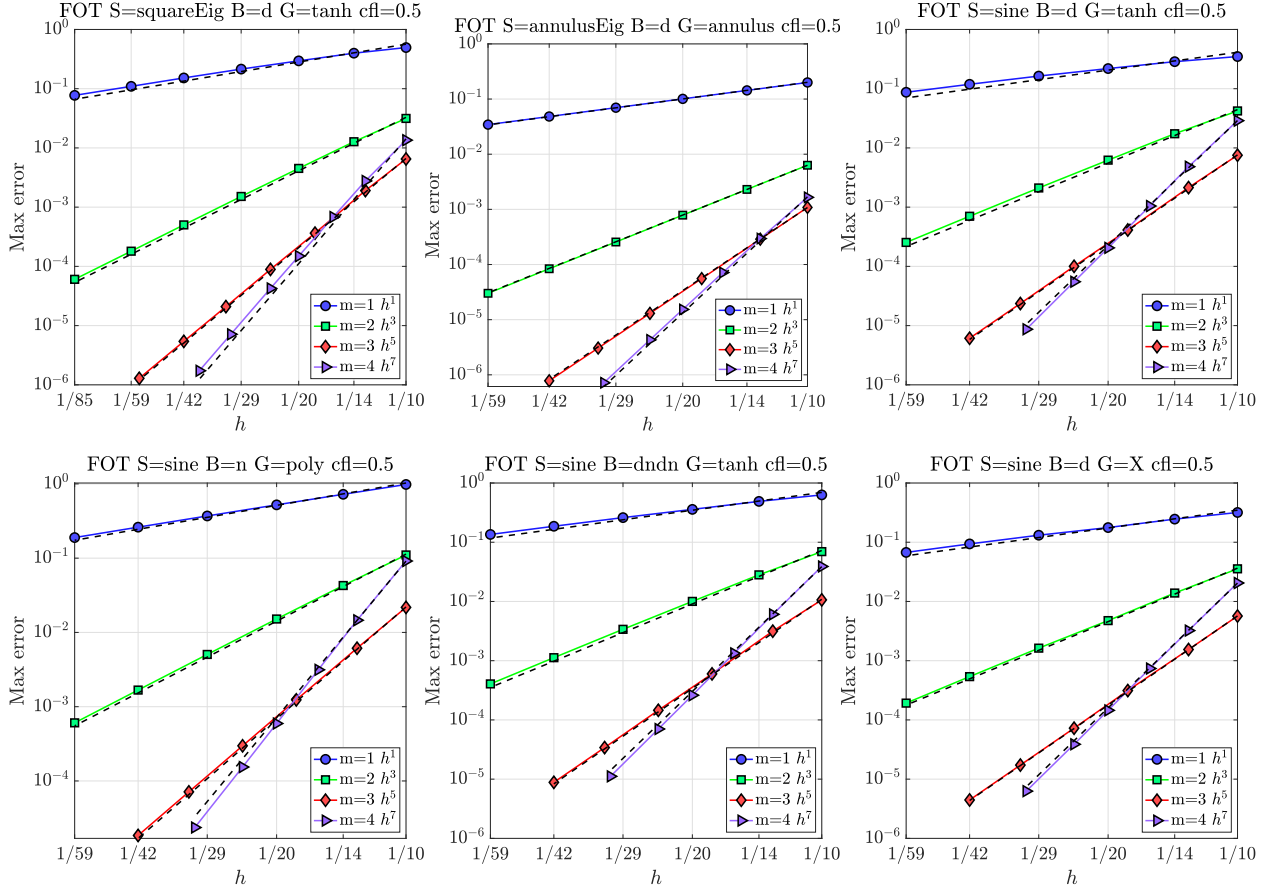
}

Grid convergence results for the FOT scheme are shown in Figure~\ref{fig:hermite2DFOT} for $m=1,2,3,4$.
These computations consider orthogonal and non-orthogonal grids, Dirichlet and Neumann boundary conditions 
(or a combination thereof to test the various treatments at corners), and both exact solutions and manufactured solutions.
In particular, results are shown for the following cases,
\begin{enumerate}
  \item Square eigenfunction, Dirichlet boundary conditions, tanh mapping,
  \item Annulus eigenfunction, Dirichlet boundary conditions, annulus mapping, 
  \item Sine solution, Neumann boundary conditions, polynomial mapping,
  \item Sine solution, Dirichlet boundary conditions, rhombus mapping,
  \item Sine solution, Dirichlet (left, bottom), Neumann (right,top), tanh mapping,
  \item Sine solution, Dirichlet boundary conditions, X mapping.
\end{enumerate}
The relative max-norm errors are computed at time $t=0.5$. The wave speed $c$ is taken as $c=1$ in all cases.
The time-step was chosen according to~\ref{sec:timeStep} with $\cfl=0.5$.
For the square eigenfunction and sine solution we take $k_x=k_y=2^{m-1}$ while for the Annulus eigenfunction
we choose the solution with $\nTheta=1+2^m-1$, and $n_r=2^{m-1}$.
In all cases the results in Figure~\ref{fig:hermite2DFOT} show that the expected order of accuracy of $2m-1$ is observed.

Figure~\ref{fig:rhombusSineFOT} shows the computed solution and errors on the rhombus grid using the FOT scheme
and the sine solution~\eqref{eq:sineSolution}.
The computed results, shown at $t=0.5$, are computed with $m=3$ using Dirichlet boundary conditions (left, bottom),
and Neumann boundary (right,top). The error is seen to be smooth up to the boundary; 
this is a good indication of accuracy and quality of the CBC conditions.

\subsection{Results for the second-order in time (SOT) scheme} \label{sec:sotResults}

{
\newcommand{\figw}{5.5cm}
\newcommand{\figh}{5.5cm}
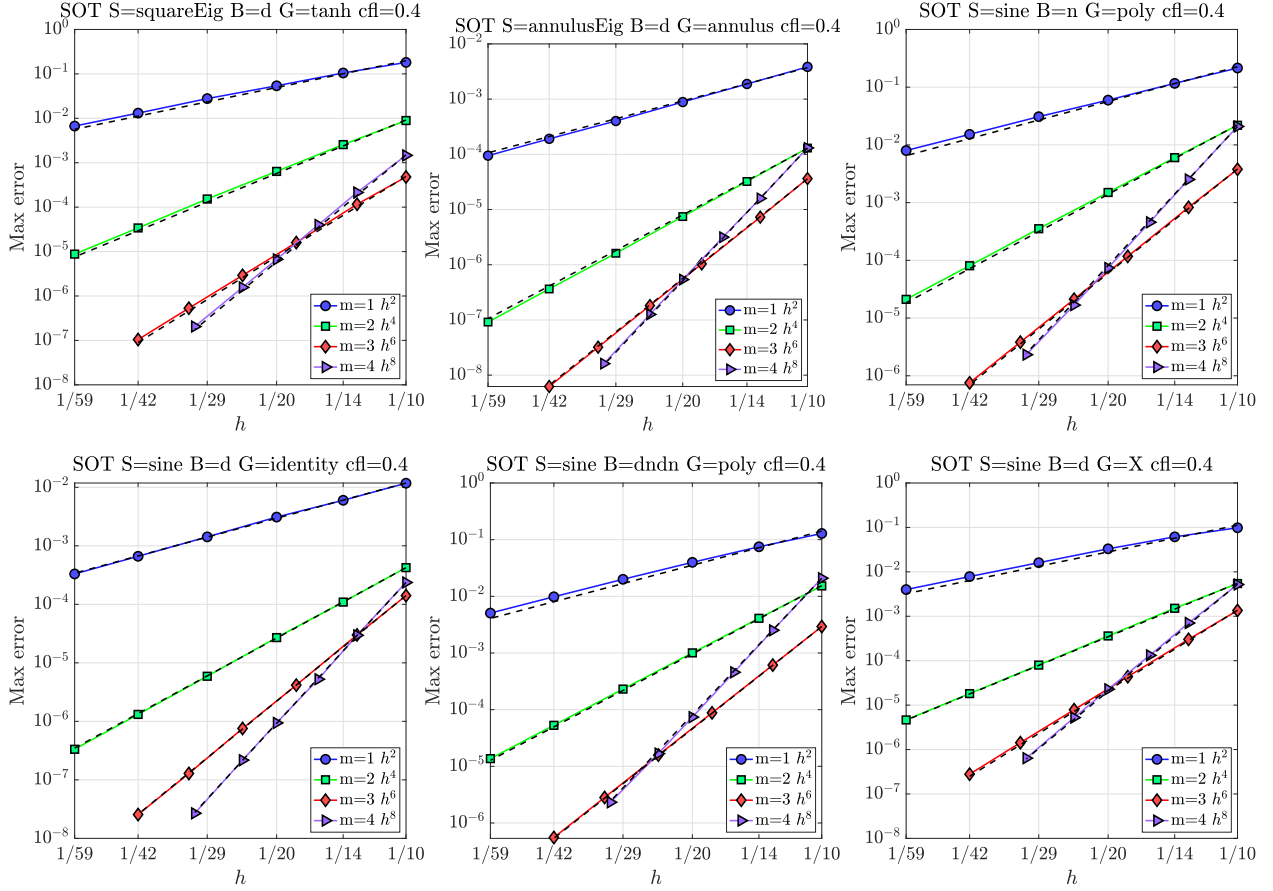
\begin{figure}[htb]
\begin{center}
\begin{tikzpicture}
  \useasboundingbox (0,.5) rectangle (16.5,12);  

  \begin{scope}[yshift=6cm]
     \figByWidth{0.0}{0}{fig/ME2DsquareEigBCdMaptConvergence}{\figw}[0][0][0][0];
     \figByWidth{5.5}{0}{fig/ME2DannulusEigBCdMapaConvergence}{\figw}[0][0][0][0];
     \figByWidth{11.}{0}{fig/ME2DsineBCnMappConvergence}{\figw}[0][0][0][0];
  \end{scope} 
  \begin{scope}[yshift=0cm]
    \figByWidth{0}{0}{fig/ME2DsineBCdMapiConvergence}{\figw}[0][0][0][0];
    \figByWidth{5.5}{0}{fig/ME2DsineBCdndnMappConvergence}{\figw}[0][0][0][0];
    \figByWidth{11.}{0}{fig/ME2DsineBCdMapxConvergence}{\figw}[0][0][0][0];
  \end{scope}   
\end{tikzpicture}
\end{center}
\caption{
Grid convergence, Hermite SOT scheme in 2D. 
Figure titles give run details: S=solution type, B=boundary conditions, G=grid type.
The expected order of accuracy for the SOT scheme is $2m$ for a degree $m$ Hermite approximation.
 }
\label{fig:hermite2DME}
\end{figure}
}

Grid convergence results for the \REV{3}{SOT} scheme are shown in Figure~\ref{fig:hermite2DFOT} for $m=1,2,3,4$.
Results are shown for the following cases,
\begin{enumerate}
  \item Square eigenfunction, Dirichlet boundary conditions, tanh mapping,
  \item Annulus eigenfunction, Dirichlet boundary conditions, annulus mapping, \REV{3}{$N_s=1$ smoothing steps (algorithm~\ref{alg:hermiteSmoothing}),}
  \item Sine solution, Neumann boundary conditions, polynomial mapping,
  \item Sine solution, Dirichlet boundary conditions, identity mapping,
  \item Sine solution, Dirichlet (left, bottom), Neumann (right,top), polynomial mapping,
  \item Sine solution, Dirichlet boundary conditions, X mapping, \REV{3}{$N_s=1$ smoothing steps (algorithm~\ref{alg:hermiteSmoothing})}.
\end{enumerate}
Parameters are chosen as in section~\ref{sec:fotResults} for the FOT scheme except that here we take $\cfl=0.4$.
In all cases the expected order of accuracy of $2m$ is observed.

Figure~\ref{fig:squareEigGridX} shows the computation of a square eigenfunction~\eqref{eq:sqaureEig} on the non-orthogonal 
X grid using the \REV{3}{SOT} scheme using $m=4$ (eight-order accurate scheme).
Figure~\ref{fig:annulusEigME} shows the computation on an eigenfunction of an annulus~\eqref{eq:annulusEig} using the 
\REV{3}{SOT} scheme with $m=5$ (tenth order scheme).
In both cases the errors are seen to be smooth up to the boundary.

\subsection{Long time simulations} \label{sec:longTimeSimulations}

\REV{3}{
To numerically demonstrate the stability of the FOT and SOT schemes some long time simulations are performed.
Results are shown for the following cases,
\begin{enumerate}
  \item Sine solution, Dirichlet (left, bottom), Neumann (right,top), identity mapping, $N_1=N_2=40$ grid points.
  \item Sine solution, Dirichlet (inner radius), Neumann (outer radius), periodic in $\theta$, annulus mapping, 
  $N_s=1$ smoothing steps (algorithm~\ref{alg:hermiteSmoothing}), $N_1=7$, $N_2=42$ grid points.
  \item Sine solution, Dirichlet (left, bottom), Neumann (right,top), X mapping, $N_s=1$ smoothing steps (algorithm~\ref{alg:hermiteSmoothing}), $N_1=N_2=20$ grid points.
\end{enumerate}
These cases test orthogonal and non-orthogonal mappings, Dirichlet and Neumann boundary conditions and all variations of corner conditions.


{
\newcommand{\figw}{6.5cm}
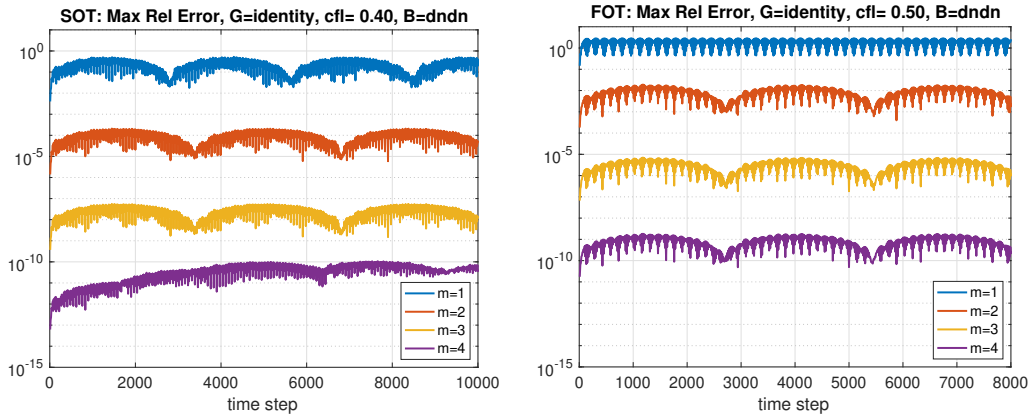
\begin{figure}[htb]
\begin{center}
\begin{tikzpicture}
  \useasboundingbox (0,.6) rectangle (14,5.5);  

  \begin{scope}[yshift=0cm]
    \figByWidth{0}{0}{fig/hermite2dStatsSchemeMEN40x40MSsineMapidentityBCdndnErrorsVersusStep}{\figw}[0.0][0][0][0];
    \figByWidth{7}{0}{fig/hermite2dStatsSchemeFOTN40x40MSsineMapidentityBCdndnErrorsVersusStep}{\figw}[0.0][0][0][0];
    
  \end{scope}   
\end{tikzpicture}
\end{center}
\caption{
Max-norm errors versus time step for the identity mapping. Left: SOT scheme. Right: FOT scheme.
 }
\label{fig:maxErrorsVersusTimeStepSquare}
\end{figure}
}

{
\newcommand{\figw}{6.5cm}
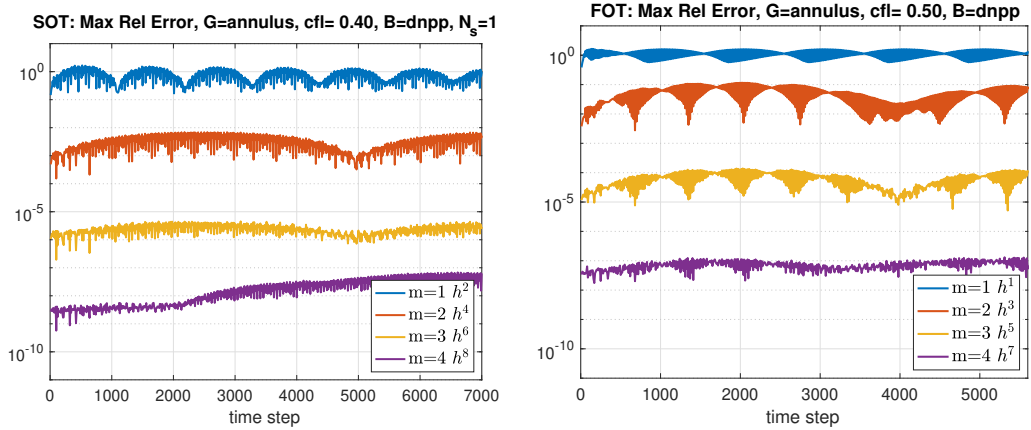
\begin{figure}[htb]
\begin{center}
\begin{tikzpicture}
  \useasboundingbox (0,.6) rectangle (14,5.5);  

  \begin{scope}[yshift=0cm]
    \figByWidth{0}{0}{fig/hermite2dStatsSchemeMEN7x7MSsineMapannulusBCdnppErrorsVersusStep}{\figw}[0.0][0][0][0];    
    \figByWidth{7}{0}{fig/hermite2dStatsSchemeFOTN7x7MSsineMapannulusBCdnppErrorsVersusStep}{\figw}[0.0][0][0][0];
    
  \end{scope}   
\end{tikzpicture}
\end{center}
\caption{
Max-norm errors versus time step for the Annulus mapping.  Left: SOT scheme. Right: FOT scheme.
 }
\label{fig:maxErrorsVersusTimeStepAnnulus}
\end{figure}
}

{
\newcommand{\figw}{6.5cm}
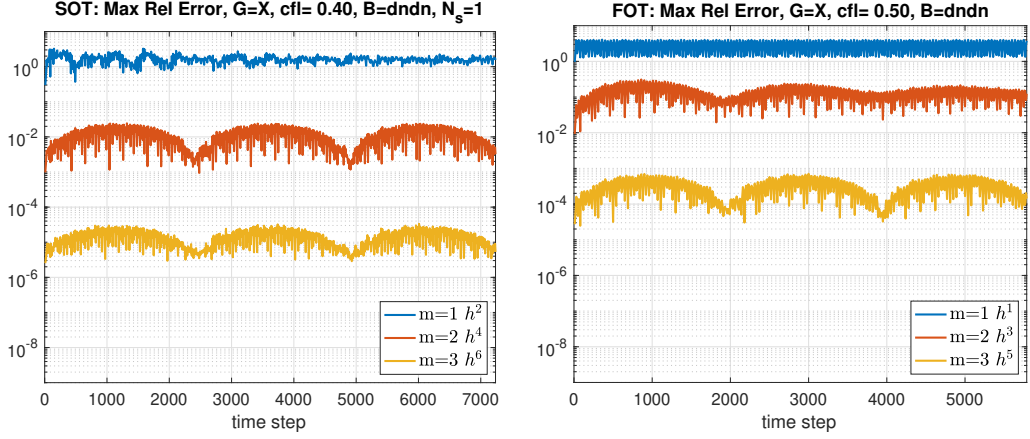
\begin{figure}[htb]
\begin{center}
\begin{tikzpicture}
  \useasboundingbox (0,.6) rectangle (14,5.5);  

  \begin{scope}[yshift=0cm]
    \figByWidth{0}{0}{fig/hermite2dStatsSchemeMEN20x20MSsineMapXBCdndnErrorsVersusStep}{\figw}[0.0][0][0][0];    
    \figByWidth{7}{0}{fig/hermite2dStatsSchemeFOTN20x20MSsineMapXBCdndnErrorsVersusStep}{\figw}[0.0][0][0][0];
    
  \end{scope}   
\end{tikzpicture}
\end{center}
\caption{
Max-norm errors versus time step for the X mapping. Left: SOT scheme. Right: FOT scheme.
 }
\label{fig:maxErrorsVersusTimeStepX}
\end{figure}
}
  
Figures~\ref{fig:maxErrorsVersusTimeStepSquare}, \ref{fig:maxErrorsVersusTimeStepAnnulus}, and \ref{fig:maxErrorsVersusTimeStepX}
show the relative max-norm errors in the solution as a function of time-step for the three cases given above.
The relative errors are computed relative to the max-norm of the true solution which is one for the sine solution.
The errors are plotted every $10$ time steps.
In each case the errors are seen to remain bounded with no exponential blowup that would be a sign of an instability.

For the X mapping results are not shown for $m=4$ since for this case large errors were detected on the Neumann boundaries
where the grid was most skew. This is consistent with the analysis in Section~\ref{sec:solvabilityCurvilinear}, 
and the determinant condition~\eqref{eq:NeumannCBCdet4} for $m=4$ which indicates that this CBC condition can be poorly conditioned on a skewed grid
when the mesh spacing is not small enough. Indeed, the errors on the Neumann boundaries for $m=4$ do decrease as the mesh is refined.

}


\section{Conclusions}  \label{sec:conclusions}

High-order accurate Hermite schemes for the wave equation on curvilinear grids have been presented. The first-order in
time (FOT) schemes have accuracy $2 m-1$ for degree $m$ Hermite schemes, while the second-order in time (SOT) 
schemes have accuracy $2 m$. Compatibility boundary conditions (CBCs) are used to build centered
polynomial approximations on the boundary. The automatic construction of the schemes for Dirichlet and Neumann boundary
conditions at arbitrary order of accuracy were given. Similarly the construction of the CBC schemes at corners was
provided. The solvability and conditioning of the matrices resulting from the CBC approximations were studied. 
\REV{4}{
With the current formulation, the conditioning of the matrices associated with the CBCs is manageable for
$m$ up to about $5$ corresponding to an order of accuracy of $10$ (SOT) or $9$ (FOT). Iterative refinement could be used to go to larger values of $m$.
It may also be possible to reformulate the equations to reduce the 
condition number. 
}
For
Cartesian grids with homogeneous Dirichlet or Neumann boundary conditions it was shown that the CBCs give polynomial
approximations with odd or even symmetry, respectively, and thus are equivalent to applying odd or even reflection
boundary conditions. Numerical examples in two dimensions demonstrated the accuracy and stability of the schemes at
different orders of accuracy and for a variety of grids, both orthogonal and non-orthogonal. Some practical
considerations in implementing the Hermite schemes on curvilinear grids were provided in the Appendices.

There are several avenues of research for future work.
The schemes \REV{3}{can} be extended to three space dimensions. 
 \REV{3}{
Choosing the CBCs at vertices, where three faces meet, will need to be worked out and the conditionning of the CBC matrices in three-dimensions is an open question.
 } 
 The
compatibility approach can be extended to treat interfaces between different material domains, where, for example the wave
speeds jump. We will investigate approaches to recover the CFL-one time-step restriction on curvilinear grids, such as
through the use of artificial dissipation or filters. Improving the conditioning of the CBC matrices for very
high-order accurate schemes could be useful. The extension of general order Hermite schemes to unstructured grids will be considered
as well as the use of the CBC approach to finite element methods.
\REV{3}{The method as presented can only handle domains that can be mapped to the unit square. Handling more complex domains would require using an overset grid or block-structured approach. This would be an an interesting
future endeavor.}

\appendix

\section{Solvability and conditioning of the CBC matrices on Cartesian grids} \label{sec:solvabilityCartesianProof}

In this section a proof of Theorem~\ref{theorem:solvabilityCartesian} is given.
The CBC matrices $M$ can be formed using a symbolic software package such as Maple.
In the case of Cartesian grids, and $m$ not too large, explicit forms of the
max-norm condition numbers can be found.
Numerical values of the condition numbers on Cartesian grids for the tall-cell ratio $\tcr=\dx/\dy=1$ can be found in Table~\ref{tab:conditionNumbersCartesian}.

\mni
\textbf{Dirichlet boundary.}
The Dirichlet CBCs~\eqref{eq:DirichletHermiteCBCs} can be row-scaled by a factor 
\ba
   C(\alpha,q) = K \, \dy^\alpha \, \dx^{2 q}, 
\ea
where $K$ is chosen to make the maximum entry (in absolute value) in the row to be one.
The resulting row-scaled equations will only depend on $\tcr=\dx/\dy$.
For $m=1$ we have an explicit formula for $\kappa_\infty(M)$ ,
\bat
   \kappa_\infty(M) & = \max(41,28+3 \tcr^2)  \times \max(\f{121}{16},1+\tcr^2) .
\eat
Note that $\kappa_\infty(M)$ only depends on the tall-cell ratio $\tcr$. As is known from other computations, the CBCs become less 
well conditioned for large $\tcr$. 
For $m=1$ and $\tcr \le 1$, the condition number is very modest in size,
\bat
   & \kappa_\infty(M) = \f{4961}{16} \approx 310 . 
\eat
For $m=2$ 
\bat
   \kappa_\infty(M) & = 
     \max \! \left(\frac{1819}{4}, 371+\frac{25 {\tcr}^{2}}{2}, \frac{597}{2}+\frac{65 {\tcr}^{2}}{4}, 5+5 {\tcr}^{4}+10 {\tcr}^{2}, \frac{1917}{16}+\frac{21 {\tcr}^{2}}{2}+\frac{15 {\tcr}^{4}}{4}\right) \\
  & \times \max \! \left(\frac{3249}{256}, 1+{\tcr}^{2}, 1+\frac{{\tcr}^{2}}{3}+{\tcr}^{4}\right) . \nonumber
\eat
For $m=2$ and $\tcr=1$ the condition number (with row-scaling) is 
\bat
   & \kappa_\infty(M) = {\frac{5909931}{1024}} \approx 5770 . 
\eat
For $m>2$ the expressions become quite lengthly and instead we just report the condition numbers for $\tcr=1$ in Table~\ref{tab:conditionNumbersCartesian}.
As $m$ increases the conditions numbers increase quite rapidly.
Even with row-scaling, the condition number for large $m$ is becoming quite large.

\mni
\textbf{Neumann boundary.} The CBCs for Neumann boundaries can also be scaled so that the resulting equations only depend on $\tcr$.
For $m=1$, 
\ba
  \kappa_\infty(M) = \max \! \left(\frac{125}{16}, 1+{ \tcr }^{2}, \frac{21}{4}+\frac{3 { \tcr }^{2}}{4}\right) \times \max \! \left(\frac{121}{16}, 1+{ \tcr }^{2}\right)
\ea
and for $\tcr\le 1$,
\ba
   \kappa_\infty(M) = \frac{15125}{256}  \approx 310 . 
\ea 
For $m=2$, 
\bat
   \kappa_\infty(M) &=
  \max \! \left(\frac{883}{16}, \frac{373}{8}+\frac{35 {\tcr}^{2}}{8}, 1+{\tcr}^{4}+2 {\tcr}^{2}, \frac{471}{32}+\frac{21 {\tcr}^{2}}{8}+\frac{15 {\tcr}^{4}}{16}\right) \\
 & \times \max \! \left(\frac{3249}{256}, 1+{\tcr}^{2}, 1+\frac{3 {\tcr}^{2}}{5}+{\tcr}^{4}\right)   \nonumber
\eat
and for $\tcr\le 1$,
\ba
   \kappa_\infty(M)  \approx 700.
\ea 
For $m>2$ and $\tcr \le 1$, see the values in Table~\ref{tab:conditionNumbersCartesian}.
It is seen that the condition numbers of the matrices for Neumann boundary conditions is similar to the condition numbers for Dirichlet boundary conditions.

\mni
\textbf{Dirichlet-Dirichlet corner.} For $m=1$, 
\ba
   \kappa_\infty(M) &= 
    \max \! \left(71, \frac{95}{2}+\frac{3}{{\tcr}^{2}}, \frac{95}{2}+3 {\tcr}^{2}, \frac{485}{16}+\frac{3 {\tcr}^{2}}{4}+\frac{3}{4 {\tcr}^{2}}\right) 
    \times \max \! \left(\frac{121}{16}, \frac{1}{3}+\frac{1}{{\tcr}^{2}}, \frac{1}{3}+{\tcr}^{2}\right) ,
\ea
and for $\tcr=1$,
\ba
   \kappa_\infty(M)  \approx 537.
\ea
The condition number at a corner grows with $\tcr^2$ and with $\tcr^{-2}$. Thus it is advisable to have $\tcr\approx 1$ near a corner.

\mni
\textbf{Neumann-Neumann corner.} For $m=1$, 
\ba
   \kappa_\infty(M) &= 
     \max \! \left(\frac{121}{16}, 1+\frac{1}{3 {\tcr}^{2}}, 1+\frac{{\tcr}^{2}}{3}, \frac{247}{64}+\frac{1}{4 {\tcr}^{2}}, \frac{247}{64}+\frac{{\tcr}^{2}}{4}, \frac{585}{256}+\frac{{\tcr}^{2}}{48}+\frac{1}{48 {\tcr}^{2}}\right) \\
    & \times \max \! \left(\frac{121}{16}, 1+\frac{1}{3 {\tcr}^{2}}, 1+\frac{{\tcr}^{2}}{3}\right), \nonumber
\ea
and for $\tcr=1$,
\ba
   \kappa_\infty(M)  \approx 572.
\ea

\mni
\textbf{Dirichlet-Neumann corner.} For $m=1$, 
\ba
   \kappa_\infty(M) &= 
     \max \! \left(\frac{43}{2}, 1+\frac{1}{{\tcr}^{2}}, \frac{115}{8}+{\tcr}^{2}, \frac{227}{16}+\frac{3}{4 {\tcr}^{2}}, \frac{573}{64}+\frac{{\tcr}^{2}}{4}+\frac{1}{16 {\tcr}^{2}}\right) \\
    & \times \max \! \left(\frac{121}{16}, 1+\frac{1}{{\tcr}^{2}}, 1+{\tcr}^{2}\right) , \nonumber
\ea
and for $\tcr=1$,
\ba
   \kappa_\infty(M)  \approx 163.
\ea

\section{Hermite evolution operators for the FOT and SOT schemes} \label{sec:evolutionAlgorithms}

The basic steps in the Hermite scheme are given in Algorithm~\ref{alg:hermite}.
The scheme involves an Hermite interpolant $\Ih$, an evolution operator $\Th$ and a boundary
condition operator $\Bh$.
In this section the evolution operators $\Th$ for the FOT and SOT schemes are described.


\subsection{FOT evolution} \label{sec:evloveFOT} 

The FOT scheme stores both the solution $u$ and the velocity $v$ on the primal and dual grids.
The degree of Taylor polynomial for $v$ is taken as one less than that for $u$.
Thus, the degree of $u_\iv(\rv)$ is $m$ and  $2m+1$ for $\uBar_\iv(\rv)$,  while the degree 
for $v_\iv(\rv)$ is $m-1$ and $2m-1$ for $\vBar_\iv(\rv)$.
Consider the process of evolving the solution on the dual grid (line 7 in Algorithm~\ref{alg:hermite})
or on the primal grid (line 10 in Algorithm~\ref{alg:hermite}).
Given $\uBar_\iv^n$ and $\vBar_\iv^{n}$, the goal is to determine $u_\iv^{n+\half}$ and  $v_\iv^{n+\half}$.
To this end, the solution is expanded in a Taylor polynomial in space and time
\bse
\label{eq:TaylorPolyFOT}
\ba
  &  \uBar_\iv^n(\rv,t) = \sum_{l_1=0}^{2m+1} \sum_{l_2=0}^{2m+1} \sum_{\ts=0}^{2m +1} \uBar_{\iv, l_1,l_2,\ts} \, R_i^{l_1} \, S_j^{l_2} \, T_n^{\ts} ,  \label{eq:TaylorPolyFOTu} \\
  &  \vBar_\iv^n(\rv,t) = \sum_{l_1=0}^{2m-1} \sum_{l_2=0}^{2m-1} \sum_{\ts=0}^{2m +1} \vBar_{\iv, l_1,l_2,\ts} \, R_i^{l_1} \, S_j^{l_2} \, T_n^{\ts} ,  \label{eq:TaylorPolyFOTv} \\
  &  R_i \eqdef \f{r - r_i}{\dr} , \quad S_j\eqdef \f{s-s_j}{\ds} ,
     \quad T_n \eqdef \f{t-t^n}{\dt} ,
\ea
\ese
for some coefficients $\uBar_{\iv, l_1,l_2,\ts}$ and $\vBar_{\iv, l_1,l_2,\ts}$.
The evolution equations for the FOT scheme enforce the following constraints 
\bse
\label{eq:FOTevolvev2d}
\bat
  & \p_x^{\alpha_1} \p_y^{\alpha_2} \p_t^\beta \p_t u = \p_x^{\alpha_1} \p_y^{\alpha_2} \p_t^\beta v    ,  
          \quad&& \alpha_k=0,1,\ldots,2 m+1,  ~~\beta=0,1,\ldots, 2m +1 , \label{eq:FOTevolveu2d} \\
  & \p_x^{\alpha_1} \p_y^{\alpha_2} \p_t^\beta \p_t v = \p_x^{\alpha_1} \p_y^{\alpha_2} \p_t^\beta( Lu  ),  
          \quad&& \alpha_k=0,1,\ldots,2 m-1, ~~\beta=0,1,\ldots,2m+1  ,
\eat
\ese
at $r=r_i$, $s=s_j$, and $t=t^n$.
The values of the coefficients $\uBar_{\iv, l_1,l_2,\ts}$ and $\vBar_{\iv,l_1,l_2,\ts}$ in~\eqref{eq:TaylorPolyFOT} for $s=0$
are determined from the Hermite interpolants for $\uBar_\iv^n$ and $\vBar_\iv^n$ at the current time $t^n$.
The values of the coefficients or $s=1,2,\ldots$ are found by a recursion derived by enforcing the conditions in~\eqref{eq:FOTevolvev2d}.
Given the coefficients  $\uBar_{\iv, l_1,l_2,\ts}$ and $\vBar_{\iv,l_1,l_2,\ts}$, the coefficients in the solution at time $t^{n+\half}$ are found
by evaluating~\eqref{eq:TaylorPolyFOT} at time $t^n+\dt/2$.
This leads to the evolution function given in Algorithm~\ref{alg:evolveFOT}.
Recall that $\Lh$ in Algorithm~\ref{alg:evolveFOT} is the matrix representation of the operator $L$. 
In the case of a Cartesian grid, enforcing~\eqref{eq:FOTevolvev2d} leads to the recursions
\bse
\ba
  & \f{\ts+1}{\dt} \, \uBar_{\iv,l_1,l_2,\ts+1} = \vBar_{\iv,l_1,l_2,\ts}, \\
  & \f{\ts+1}{\dt} \, \vBar_{\iv,l_1,l_2,\ts+1} = c^2 \f{(l_1+2)(l_1+1)}{\dx^{2}} \, \uBar_{\iv,l_1+2,l_2,\ts} + c^2 \f{(l_2+2)(l_2+1)}{\dy^{2}}\,  \uBar_{\iv,l_1,l_2+2,\ts} ,
\ea
\ese
for $\ts=0,1,2,\ldots,2m+1$.
{
\newcommand{\algFontSize}{\small}
\begin{algorithm}[hbt]
\algFontSize 
\caption{FOT Evolution: advance the solution for half a time-step.}
\begin{algorithmic}[1]
  \Function{$[u_\iv^{n+\half}, v_\iv^{n+1/2}] = $evolveFOT}{$\uBar_\iv^{n}, \vBar_\iv^{n}$}
    \State $\uBar_{\iv,l_1,l_2,0} = \uBar_{\iv,l_1,l_2}^{n}, \qquad l_1,l_2=0,1,2,\ldots,2 m+1$
    \State $\vBar_{\iv,l_1,l_2,0} = \vBar_{\iv,l_1,l_2}^{n}, \qquad~ l_1,l_2=0,1,2,\ldots,2 m-1$
    \For{$\ts=0,1,\ldots,2 m +1$} 
        \State $\displaystyle \uBar_{\iv,l_1,l_2,\ts+1} = \f{\dt}{s+1} \vBar_{\iv,l_1,l_2,\ts} , \qquad l_1,l_2=0,1,2,\ldots,2 m+1$
        \State $\wBar_{\iv,l_1,l_2} =  \uBar_{\iv,l_1,l_2,\ts} ,                         \hspace{58pt} l_1,l_2=0,1,2,\ldots,2 m+1 $                                      
        \State $\displaystyle \vBar_{\iv,l_1,l_2,\ts+1} = \f{\dt}{s+1} (\Lh \wBar_{\iv})_{\iv,l_1,l_2,\ts}$ 
           \Comment See Algorithm~\ref{alg:applyOperator} for $\Lh \wBar_{\iv}$
    \EndFor 
    \State $\displaystyle u_\iv^{n+\half,l_1,l_2}   = \sum_{\ts=0}^{2m +2} \uBar_{\iv, l_1,l_2,\ts} \, \Big(\half\Big)^{\ts} , \qquad l_1,l_2=0,1,2,\ldots,m+1$
     \Comment Evaluate Taylor series in time
    \State $\displaystyle v^{n+\half}_{\iv,l_1,l_2} = \sum_{\ts=0}^{2m +2} \vBar_{\iv, l_1,l_2,\ts} \, \Big(\half\Big)^{\ts} , \qquad \hspace{16pt} l_1,l_2=0,1,2,\ldots,m-1 $
    \Comment Evaluate Taylor series in time
 \EndFunction
\end{algorithmic} 
\label{alg:evolveFOT}
\end{algorithm}
}

\subsection{SOT evolution} \label{sec:evolveME}

The evolution of the SOT scheme is based on the Taylor series expansion of the second divided difference in time,
\ba
   \f{u(\rv,t+\delta)-2u(\rv,t)+u(\rv,t-\delta)}{\delta^2} = 2 \sum_{\mu=1}^\infty \f{1}{(2\mu)!} \delta^{2\mu}  \p_t^{2\mu} u(\rv,t).
\ea
Using $\p_t^{2\mu} u = L^\mu \mu$ and setting $\delta=\dt/2$ leads to
\ba
  u(\rv,t+\f{\dt}{2})= 2 u(\rv,t) -  u(\rv,t-\f{\dt}{2}) + 2 \sum_{\mu=1}^\infty \f{1}{(2\mu)!} \left[\f{\dt}{2}\right]^{2\mu} L^\mu u(\rv,t). \label{eq:TaylorME}
\ea
The SOT evolution equations are derived from taking spatial derivatives of~\eqref{eq:TaylorME}, leading to the approximations
\ba
  & \p_{r_1}^{\alpha_1}\p_{r_2}^{\alpha_2} \, u(\rv,t^{n+\half})  
  \approx  \nonumber \\
  &\qquad  2\,  \p_{r_1}^{\alpha_1}\p_{r_2}^{\alpha_2} \, u(\rv,t^n) -  \p_{r_1}^{\alpha_1}\p_{r_2}^{\alpha_2} \, u(\rv,t^{n-\half})
         + 2 \sum_{\mu=1}^{m}  \f{1}{(2\mu)!} \left[\f{\dt}{2}\right]^{2\mu}  \p_{r_1}^{\alpha_1}\p_{r_2}^{\alpha_2}\, L^{\mu} u(\rv,t),
\label{eq:evolveME}
\ea 
at $\rv=\rv_\iv$ and for $\alpha_1,\alpha_2=0,1,2,\ldots,m$.
Algorithm~\ref{alg:evolveME} gives the SOT evolution function used to compute the Hermite DOFs $u_{\iv,l_1,l_2}^{n+\half}$ based on~\eqref{eq:evolveME}.
{
\newcommand{\algFontSize}{\small}
\begin{algorithm}
\algFontSize 
\caption{SOT Evolution: advance the solution a half time-step.}
\begin{algorithmic}[1]
  \Function{$ u_\iv^{n+\half} = $ evolveSOT}{$\uBar_\iv^n$, $u_\iv^{n-\half}$}
    \State $u_{\iv,l_1,l_2}^{n+\half} = 2 \, \uBar_{\iv,l_1,l_2}^{n} - u_{\iv,l_1,l_2}^{n-\half}, \qquad l_1,l_2=0,1,2,\ldots,m$
    \State $\wBar_\iv  = \uBar_\iv^{n}$  \Comment Holds $L^\mu \uBar_\iv^{n}$
    \For{$\mu=1,2,\ldots,m$} 
        \State $\wBar_\iv = \Lh \wBar_\iv$                               \Comment See Algorithm~\ref{alg:applyOperator} for $\Lh \wBar_{\iv}$    
        \State $\displaystyle u_{\iv,l_1,l_2}^{n+\half} = u_{\iv,l_1,l_2}^{n+\half} + 2 \f{1}{(2\mu)!} \left[\f{\dt}{2}\right]^{2\mu} \wBar_{\iv,l_1,l_2}, \qquad l_1,l_2=0,1,2,\ldots,m$
    \EndFor 
 \EndFunction
\end{algorithmic} 
\label{alg:evolveME}
\end{algorithm}
}

\mni
The explicit form of the update on a Cartesian grid is 
\shadedBoxWithShadow{align}{orange}{
&    u_{\iv,l_1,l_2}^{n+\half} 
   = 2  \,\uBar_{\iv,l_1,l_2}^{n} - u_{\iv,l_1,l_2}^{n-\half} \nonumber \\
 &\quad    + \f{2}{l_1! \, l_2!} \sum_{\mu=1} \f{1}{(2 \mu)!} \Big[\f{c \dt}{2}\Big]^{2 \mu} 
      \sum_{j=0}^\mu  { \mu \choose j} \f{(2(\mu-j)+l_1)!}{h_x^{2(\mu-j)}} \, \f{(2j+l_2)!}{h_y^{2j}}
          \uBar_{\iv,2(\mu-j)+l_1, 2j+l_2}^n, 
}
for $l_1,l_2=0,1,\ldots,m$.

\subsubsection{First time-step}

  The SOT scheme requires two starting values, the solution on the primal grid at $t=0$ and 
the solution on the dual grid at $t=-\dt/2$. The solution and it's spatial derivatives at $t=0$ are found 
from the initial condition~\eqref{eq:waveU0}.
The solution at $t=-\dt/2\eqdef \delta$ can be found from a Taylor series in time,
\ba
   u(\xv,t+\delta) = u(\xv,0) + \delta \, \p_t u(\xv,0) + \f{\delta^2}{2!} \,\p_t^2 u(\xv,0)
          + \f{\delta^3}{3!}\, \p_t^3 u(\xv,0)+ \f{\delta^4}{4!}\, \p_t^4 u(\xv,0) + \ldots  \label{eq:initialStepTaylor}
\ea
Given initial conditions,
\bse
\ba
  &  u(\xv,0)      = U_0(\xv), \\
  &  \p_t u(\xv,0) = U_1(\xv), 
\ea
\ese
the even time-derivatives of $u(\xv,0)$ are (assuming here that the body forcing $f(\xv,t)$ is zero)
\ba
   \p_t^{2q} u(\xv,0) = L^q U_0(\xv), \qquad q=0,1,2,\ldots ,
\ea
where $L=c^2\Delta$. The odd time-derivatives are 
\ba
   \p_t^{2q+1} u(\xv,0) = L^q U_1(\xv) , \qquad q=0,1,2,\ldots .
\ea
The Taylor series~\eqref{eq:initialStepTaylor} can also be used to update the spatial derivatives.
Algorithm~\ref{alg:takeFirstStep} outlines the first (backward) step.
For degree $m$ one should keep $2m +1$ terms in the Taylor series (the last term is the $\delta^{2m}$ term).
The Algorithm~\ref{alg:takeFirstStep} keeps one additional term.

{
\newcommand{\algFontSize}{\small}
\begin{algorithm}
\algFontSize 
\caption{First (backward) time-step for the SOT scheme}
\begin{algorithmic}[1]
  \Function{[$u_\iv^0$, $u^{-\half}_\jv$] = takeFirstStep}{}  
    \State Set $\displaystyle u_\iv^0$  from $\p_{r_1}^{l_1}\p_{r_2}^{l_2} U_0, \hspace{29pt} l_1,l_2=0,1,\ldots,m, \quad \iv\in P$  \Comment Initial solution.
    \State Set $\displaystyle v_\iv^0$  from $\p_{r_1}^{l_1}\p_{r_2}^{l_2} U_1, \hspace{30pt} l_1,l_2=0,1,\ldots,m, \quad \iv\in P$  \Comment Initial time-derivative.
    \State $\uBar_\jv  = \Ih( u_\iv^0 )$, \hspace{64pt} $\iv\in P$, $\jv\in D$                 \Comment Interpolate to dual grid.
    \State $\vBar_\jv  = \Ih( v_\iv^0 )$, \hspace{65pt} $\iv\in P$, $\jv\in D$                 \Comment Interpolate to dual grid.
    \State $\delta = -\dt/2$
    \State $u^{-\half}_{\jv,l_1,l_2} = \uBar_{\jv,l_1,l_2} + \delta \, \vBar_{\jv,l_1,l_2}, \quad l_1,l_2=0,1,\ldots,m, \quad \jv\in D$
    \For{$k=1,2,\ldots,m$} 
       \State $\displaystyle u^{-\half}_{\jv,l_1,l_2} = u^{-\half}_{\jv,l_1,l_2} 
           + \f{\delta^{2k}}{2k !}\, (\Lh^k \uBar_{\jv})_{\jv,l_1,l_2} + \f{\delta^{2k+1}}{(2k+1) !}\, (\Lh^k \vBar_{\jv})_{\jv,l_1,l_2},  \quad l_1,l_2=0,1,\ldots,m, \quad \jv\in D$ 
    \EndFor    
 \EndFunction

\end{algorithmic} 
\label{alg:takeFirstStep}
\end{algorithm}

}

\subsection{Choosing the time-step} \label{sec:timeStep}

On Cartesian grids, the time step is chosen from 
\ba
   \f{c \, \dt}{\min({\dx,\dy})} = \cfl,
\ea
where $\cfl$ is the CFL parameter.
The SOT schemes with CBCs on Cartesian grids appear to be stable to $\cfl=1.0$ (found experimentally).
The FOT scheme on Cartesian grids requires special fixes to reach $\cfl=1.0$ as described in~\cite{alvarez2022hermite}.
%
On curvilinear grids we estimate the smallest grid spacing in the $r$ and $s$ coordinate directions 
from the grid points, 
\ba
  \dxMin = \min_{\iv} \Big( | \xv_{i_1+1,i_2} - \xv_\iv | , | \xv_{i_1,i_2+1} - \xv_\iv | \Big),
\ea
and choose $\dt$ from 
\ba
   \f{c \, \dt}{\dxMin} = \cfl .
\ea
On curvilinear grids, the schemes generally have a lower CFL limit than for Cartesian grids 
 (but as the mesh is refined this limit appears to approach the Cartesian grid stable CFL). 
For the computations in this article a choice of $\cfl=0.5$ was taken for the FOT scheme and $\cfl=0.4$ for the SOT scheme, unless
otherwise specified.
Numerical experiments suggest that the addition of some dissipation to either scheme will increase
the stable $\cfl$. An investigation into this behavior will be left to future work.

\section{Practicalities}   \label{sec:practicalities}

This section provides some helpful information for those readers interested in implementing Hermite schemes.

\subsection{Hermite interpolants} \label{sec:HermiteInterpolants}

In one space dimension the degree $m$ Taylor polynomial representation of the solution is 
\ba
  &  u_i(r) = \sum_{l_1=0}^{m}  u_{i,l_1} \, R_i^{l_1}, \qquad  R_i \eqdef \f{r - r_i}{\dr}  .  \label{eq:TaylorPoly1D} 
\ea
The degree $2m+1$ Hermite interpolant given by
\ba
  &  \uBar_{i+\half}(r) = \sum_{l_1=0}^{2m+1}  \uBar_{i+\half,l_1} \, R_i^{l_1} , 
\ea
is chosen to match the solution and it's derivatives at points $r_i$ and $r_{i+1}$,
\bse
\bat
    & \p_r^\alpha \uBar_{i+\half}(r_i    ) = \p_r^\alpha u_i(r_i) ,     \quad&& \alpha=0,1,2,\ldots,m, \\
    & \p_r^\alpha \uBar_{i+\half}(r_{i+1}) = \p_r^\alpha u_{i+1}(r_{i+1}) , \quad&& \alpha=0,1,2,\ldots,m.
\eat
\ese
Now
\ba
   & \f{\dr^\alpha}{\alpha!} \p_r^\alpha u_i(r_i)     = u_{i,  \alpha}, \quad 
     \f{\dr^\alpha}{\alpha!} \p_r^\alpha u_{i+1}(r_{i+1}) = u_{i+1,\alpha}, 
\ea
while
\ba
   \f{\dr^\alpha}{\alpha!} \p_r^\alpha \uBar_{i+\half}(r) = \sum_{l_1=\alpha}^{2m+1} {l_1 \choose \alpha } \uBar _{i+\half,l_1} \, R_i^{l_1-\alpha}  .
\ea
This leads to the interpolation conditions
\bse
\label{eq:HermiteInterpolation1d}
\ba
  & \sum_{l_1=\alpha}^{2m+1} {l_1 \choose \alpha } \uBar _{i+\half,l_1} \, \Big[ -\half \Big]^{l_1-\alpha}  = u_{i,  \alpha}, \\
  & \sum_{l_1=\alpha}^{2m+1} {l_1 \choose \alpha } \uBar _{i+\half,l_1} \, \Big[ +\half \Big]^{l_1-\alpha}  = u_{i+1,\alpha} .
\ea
\ese
for $\alpha=0,1,2,\ldots,m$. Equations~\eqref{eq:HermiteInterpolation1d} define a linear system of equations whose solution
can be written as 
\ba
   \uBar_{i+\half} = \Ih^{(1)}( u_i ),
\ea
where $\Ih^{(1)}$ is the Hermite interpolation operator in coordinate direction $r_1$.
The Hermite interpolant in two dimensions is defined by repeated application of one-dimensional interpolants,
first in the $r_1$-direction and then in the $r_2$-direction as given in in Algorithm~\ref{alg:hermiteInterpolant2d}.
{
\newcommand{\hermiteInterpolantOneD}{{\blue hermiteInterpolantOneD}}
\newcommand{\algFontSize}{\small}
\begin{algorithm}
\algFontSize 
\caption{Compute the Hermite interpolant in two dimensions, $\uBar_\iv = \Ih( u_\iv)$.}
\begin{algorithmic}[1]
  \Function{$ \uBar_\iv = \Ih$}{ $u_\iv$ }
  \For{ $l_2=0,1,\ldots,m$ } \Comment Interpolate in $r_1$ direction
    \State  $\uBar_{i+\half,i_2, 0:2m+1,l_2 } = \Ih^{(1)}( u_{i_1,i_2, 0:m+1,l_2} ), \hspace{37pt} i_1=0,1,\ldots,N_1-1, ~~ i_2=0,1,\ldots,N_2$
  \EndFor
  \For{ $l_1=0,1,\ldots,2m+1$ } \Comment Interpolate in $r_2$ direction
     \State $\uBar_{i_1+\half,i_2+\half, l_1,0:2m+1} = \Ih^{(2)}( \uBar_{i_1+\half,i_2, l_1,0:m+1} ), \quad i_1=0,1,\ldots,N_1-1, ~~ i_2=0,1,\ldots,N_1-1 $
  \EndFor
  \EndFunction
\end{algorithmic} 
\label{alg:hermiteInterpolant2d}
\end{algorithm}
}

\subsection{Taylor polynomial coefficients from function evaluations} \label{sec:TaylorFromFunctioValues}

Algorithm~\ref{alg:getTaylorCoefficients} contains a useful procedure that 
determines approximations to the scaled Taylor polynomial coefficients (as used by the Hermite schemes) of a known function $f(x)$.
This function can be used for setting up initial conditions, evaluating forcing functions and their derivatives for boundary conditions 
(such as in \eqref{eq:DirichletHermiteCBCs}), as well as
computing the Taylor polynomial representations of the curvilinear coefficients of $L$ in ~\eqref{eq:WaveParameterSpace}.

Algorithm~\ref{alg:getTaylorCoefficients} finds approximate Taylor coefficients for a function $f(x)$ on a cell with center $x_i$ and width $\dx$.
The algorithm evaluates the function on a set of $q+1$ scaled Chebyshev points on $[x_i-\dx/2,x_i+\dx/2]$, forms an interpolant in Newton divided difference form,
and then converts the coefficients of this interpolant into the coefficients of a power series.
This procedure is described in \textsl{Algorithm for Dual System}, in Solution of Vandermonde Systems of Equations, by {\AA}ke Bj\"ork and Victor Pereyra~\cite{BjorckPereyra1970}.
It avoids computing the interpolant from a Vandermonde system, which would be very ill conditioned.
For a Hermite polynomial of degree $m$ we choose $q= 2m + 1$ interpolation intervals.

{
\newcommand{\algFontSize}{\small}
\begin{algorithm}
\algFontSize 
\caption{Find scaled Taylor coefficients $f_j$, $j=0,1,\ldots,q$ from function evaluations.}
\begin{algorithmic}[1]
  \Function{$ \fv  $ = getTaylorCoefficients}{ $q$, $x_i$, $\dx$, $f$ }
    \State // Evaluate the function $f(x)$ at $q+1$ Chebyshev points on the cell centered at $x_i$ with width $\dx$
    \For{ $j=0,1,2,\ldots,q$ }
      \State $z_j = - \half \cos(\pi j /q)$  \Comment Chebyshev points on $[-\half,\half]$
      \State $f_j = f(x_i + \dx \, z_j)$  \Comment Evaluate $f$ on scaled Chebyshev points on $[x_i-\dx/2,x_i+\dx/2]$
    \EndFor

    \For{ $k=1,2,\ldots,q$ }
      \For{ $j=q,q-1,\ldots,k$ }
        \State $\displaystyle f_j = \f{f_j - f_{j-1}}{z_j - z_{j-k}}$ \Comment Newton divided differences
      \EndFor
    \EndFor

   \For{ $k=q-1,q-2,\ldots,0$ }
     \For{ $j=k,k+1,\ldots,q-1$ }  
       \State $f_j = f_j - z_k \, f_{j+1}$  \Comment Recursion to convert Newton form to Taylor form
     \EndFor
   \EndFor   
  \EndFunction  
\end{algorithmic} 
\label{alg:getTaylorCoefficients}
\end{algorithm}
}

\subsection{Taylor polynomials for the Laplacian in curvilinear coordinates}
\label{sec:TaylorForLaplacianInCurvilinear}
The coefficients, $a^{\alpha\,\beta}(\rv)$, in the Laplacian in curvilinear coordinates in~\eqref{eq:WaveParameterSpace},
are approximated as Taylor polynomials of degree $2m+1$ in each direction,
\ba
   a^{\alpha\,\beta}_{\iv}(\rv) = \sum_{l_1=0}^{2m+1} \sum_{l_2=0}^{2m+1} a^{\alpha\,\beta}_{\iv,l_1,l_2}\, R_i^{l_1} \, S_j^{l_2}.
\ea
Values for $a^{\alpha\,\beta}_{\iv,l_1,l_2}$ can be conveniently computed making use of the approach
outlined in~\ref{sec:TaylorFromFunctioValues} which requires a function to evaluate 
$a^{\alpha\,\beta}_{\iv}(\rv)$ for different values of $\rv$.
Consider, for example, computing the Taylor polynomial approximation to $a^{2 0}(\rv)$,
\ba
  a^{20}(\rv) = (\p_{x_1} r_1)^2 + (\p_{x_2} r_1)^2 . \label{eq:a20}
\ea
Let us suppose that we have a function to compute entries in the Jacobian matrix\footnote{If the metric derivatives are not known then one can work directly with
$\Gv(\rv)$ and take
derivatives of the Taylor polynomials to get the entries in the Jacobian matrix.}
\ba
   \f{\p \xv}{\p \rv} = \f{\p\Gv(\rv)}{\p \rv} = \Big[ \f{\p x_\mu }{\p r_\nu} \Big].
\ea
The inverse metrics can be found from the inverse of the Jacobian matrix
\ba
   \f{\p \rv}{\p \xv} = \Big[\f{\p \xv}{\p \rv}\Big]^{-1} = \Big[ \f{\p r_\mu }{\p x_\nu} \Big].
\ea
Now  fit Taylor polynomials to the inverse metrics
\ba
   \Big[\f{\p r_\mu }{\p x_\nu}\Big]_{\iv}(\rv) 
      = \sum_{l_1=0}^{2m+1} \sum_{l_2=0}^{2m+1} \Big[\f{\p r_\mu }{\p x_\nu}\Big]_{\iv,l_1,l_2}\, R_i^{l_1} \, S_j^{l_2}.
\ea
The coefficients in the Taylor polynomial for $a^{20}(\rv)$ in~\eqref{eq:a20} can then be computed
from the polynomials for the inverse metrics through multiplication and addition of polynomials, truncated to degree $2m+1$
in each direction.
The coefficients $a^{10}$ and $a^{01}$ in~\eqref{eq:WaveParameterSpace} depend on derivatives of the inverse metrics and these 
can be computed using the chain rule and taking derivatives of the Taylor polynomials.

\subsection{Applying the wave operator in curvilinear coordinates} \label{sec:applyingTheWaveOperatorInCurvilinearCoordinates}

Given the Taylor polynomial representations for the curvilinear coefficients $a^{\mu,\nu}_{\iv}(\rv)$
applying the wave operator $L$ to a Hermite representation $\uBar_{\iv}(\rv)$ is straightforward.
For example, consider the computation of $\wBar = a^{20} \p_r^2 \uBar$. 
We have
\ba
   \p_r^2 \uBar_{\iv}(\rv) = \sum_{l_1=2}^{2m+1} \sum_{l_2=0}^{2m+1} \uBar_{\iv,l_1,l_2}\, \f{l_1(l_1-1)}{\dr^2} R_i^{l_1-2 } \, S_j^{l_2}.
\ea
Then
\bse
\ba
  a^{20}_{\iv}(\rv) \p_r^2 \uBar_\iv(\rv) &= 
             \sum_{l_1'=0}^{2m+1} \sum_{l_2'=0}^{2m+1} a^{2 0}_{\iv,l_1',l_2'}\, R_i^{l_1'} \, S_j^{l_2'} 
      \times \sum_{l_1''=2}^{2m+1} \sum_{l_2''=0}^{2m+1} \uBar_{\iv,l_1'',l_2''}\, \f{l_1''(l_1''-1)}{\dr^2} R_i^{l_1''-2 } \, S_j^{l_2''}, \\
      &=  \sum_{l_1'=0}^{2m+1} \sum_{l_2'=0}^{2m+1} \sum_{l_1''=2}^{2m+1} \sum_{l_2''=0}^{2m+1}
             a^{2 0}_{\iv,l_1',l_2'}\, \uBar_{\iv,l_1'',l_2''}\, \f{l_1''(l_1''-1)}{\dr^2} R_i^{l_1'+l_1''-2} \, S_j^{l_2'+l_2''}
\ea
\ese
Setting $l_1'+l_1''-2=l_1$ and $l_2'+l_2''=l_2$ and $l_1'=k_1$ and $l_2'=k_2$ implies 
\ba
    l_1''=l_1-k_1+2, \quad l_2'' = l_2-k_2
\ea
and thus 
\bse
\label{eq:a20TaylorPoly}
\ba
  &  a^{20}_{\iv}(\rv) \p_r^2 \uBar_\iv(\rv) =  \sum_{l_1=0}^{2m+1} \sum_{l_2=0}^{2m+1} \wBar_{\iv,l_1,l_2} \, R_i^{l_1} \, S_j^{l_2}, \\
  & \wBar_{\iv,l_1,l_2} = \sum_{k_1=0}^{l_1} \sum_{k_2=0}^{l_2}
         a^{\alpha\,\beta}_{\iv,k_1,k_2}\, \uBar_{\iv,l_1-k_1+2,l_2-k_2}\, \f{(l_1-k_1+2)(l_1-k_1+1)}{\dr^2}
\ea
\ese
Note that the polynomials in~\eqref{eq:a20TaylorPoly} have been truncated to degree $2m+1$ and that any terms with an array index out of bounds should be ignored.
Algorithm~\ref{alg:applyOperator} gives the full algorithm to apply $L$ (with the same caveat that terms with invalid subscripts should be ignored.)

{
\newcommand{\algFontSize}{\small}
\begin{algorithm}
\algFontSize 
\caption{Evaluate $\wBar_\iv = \Lh \uBar_\iv$, for $L$ in curvilinear coordinates. }
\begin{algorithmic}[1]
  \Function{$ \wBar_\iv  $ = applyOperator}{ $\uBar_\iv$ }
  \State // $L = c^2 \Delta = a^{20}(\rv) \frac{\p^2}{\p r_1^2} + a^{11}(\rv)\frac{\p^2}{\p r_1 \p r_2} + a^{02}(\rv) \frac{\p^2}{\p r_2^2} 
         + a^{10}(\rv)\frac{\p}{\p r_1} + a^{01}(\rv) \frac{\p }{\p r_2}$
  \For{ $l_2=0,1,\ldots,2m+1$ } 
  \For{ $l_1=0,1,\ldots,2m+1$ } 
    \State $\wBar_{\iv,l_1,l_2}=0$
    \For{ $k_2=0,1,\ldots,l_1$ } 
    \For{ $k_1=0,1,\ldots,l_2$ }  
      \State $\displaystyle \wBar_{\iv,l_1,l_2} = \wBar_{\iv,l_1,l_2} $
      \State $\displaystyle \hspace{30pt} + a^{20}_{\iv,k_1,k_2} \uBar_{\iv,l_1-k_1+2,l_2-k_2  } \f{(l_1-k_1+2)(l_1-k_1+1)}{\dr^2} $
      \State $\displaystyle \hspace{30pt} + a^{11}_{\iv,k_1,k_2} \uBar_{\iv,l_1-k_1+1,l_2-k_2+1} \f{(l_1-k_1+1)(l_2-k_2+1)}{\dr\ds}$
      \State $\displaystyle \hspace{30pt} + a^{02}_{\iv,k_1,k_2} \uBar_{\iv,l_1-k_1  ,l_2-k_2+2} \f{(l_2-k_2+2)(l_2-k_2+1)}{\ds^2} $
      \State $\displaystyle \hspace{30pt} + a^{10}_{\iv,k_1,k_2} \uBar_{\iv,l_1-k_1+1,l_2-k_2  } \f{(l_1-k_1+1)}{\dr}  + a^{01}_{\iv,k_1,k_2} \uBar_{\iv,l_1-k_1  ,l_2-k_2+1} \f{(l_2-k_2+1)}{\ds}$
    \EndFor
    \EndFor
  \EndFor
  \EndFor
 \EndFunction
\end{algorithmic} 
\label{alg:applyOperator}
\end{algorithm}
}

\bibliographystyle{elsart-num}
\bibliography{hcbc}

\end{document}